\documentclass[11pt,a4paper]{article}
\usepackage[italian,english]{babel}
\usepackage{amsfonts}
\usepackage{amsmath}
\usepackage{amsthm}
\usepackage{amssymb}
\usepackage{enumerate}
\usepackage{xspace}
\usepackage{euscript}
\usepackage{graphicx}
\usepackage{amscd}
\usepackage{epsfig}
\usepackage{tabularx}
\usepackage{diagrams}
\usepackage[all]{xy}

        \newtheorem{lemma}{Lemma}[section]
        \newtheorem{proposition}[lemma]{Proposition}
        \newtheorem{theorem}[lemma]{Theorem}
        \newtheorem{corollary}[lemma]{Corollary}
        \newtheorem{definition}[lemma]{Definition}
        
      \theoremstyle{definition}

        \newtheorem{rem}[lemma]{Remark}

\begin{document}

\title{\textbf{Quantum Cohomology of Hilb$^2(\mathbb{P}^1\times\mathbb{P}^1)$
    and Enumerative Applications}}
\author{Dalide Pontoni}
\date{ }
\maketitle
\addtolength{\topmargin}{-1.0cm}
\addtolength{\textheight}{2.0cm}

\footskip=40pt

\newcommand{\g}{\ensuremath{\mathbf{G}}}
\newcommand{\p}[1]{\ensuremath{\mathbb{P}(#1)}}
\newcommand{\ag}[1]{\ensuremath{A^{#1}(\g)}}
\newcommand{\e}{\ensuremath{\mathcal{E}}}
\newcommand{\h}{\ensuremath{\mathbf{H}}}
\newcommand{\s}{\ensuremath{\mathbf{S}}}
\newcommand{\D}{\ensuremath{\Delta}}
\newcommand{\ah}[1]{\ensuremath{A^{#1}(\h)}}
\newcommand{\sig}[1]{\ensuremath{\sigma_{#1}}}               
\newcommand{\w}[1]{\ensuremath{W_{#1}}}
\newcommand{\wt}[1]{\ensuremath{\tilde{W}_{#1}}}
\newcommand{\aq}[1]{\ensuremath{A^{#1}(Q)}}
\newcommand{\G}[1]{\ensuremath{\Gamma(#1)}}
\newcommand{\gm}[1]{\ensuremath{\Gamma_{#1}}}
\newcommand{\mo}{\ensuremath{\Sigma}}
\newcommand{\omo}{\ensuremath{\overline{\mo}}}
\newcommand{\td}{\ensuremath{\tilde{\delta}}}
\newcommand{\tu}{\ensuremath{\tilde{U}}}
\newcommand{\cuno}{\ensuremath{C_1}}
\newcommand{\cdue}{\ensuremath{C_2}}
\newcommand{\ctre}{\ensuremath{F}}
\newcommand{\Ll}{\ensuremath{\Lambda(l)}}
\newcommand{\Lp}{\ensuremath{\Lambda(p)}}
\newcommand{\supp}[1]{\ensuremath{\textrm{Supp}~{#1}}}
\newcommand{\mbar}[1]{\ensuremath{\overline{M}_{#1}}}
\newcommand{\ed}[1]{\ensuremath{\textrm{ed}_{\scriptscriptstyle{#1}}}}
\newcommand{\ds}[1]{\ensuremath{\textrm{d}_{\scriptscriptstyle{#1}}}}
\newcommand{\ext}{\ensuremath{\textrm{Ext}}}
\newcommand{\R}{\ensuremath{\textrm{R}}}
\newcommand{\cod}{\ensuremath{\textrm{cod}~}}
\newcommand{\dg}{\ensuremath{\textrm{deg}~}}
\newcommand{\m}[1]{\ensuremath{\mathcal{M}_{#1}}}

\abstract{
\noindent
We compute the Small Quantum Cohomology of
$\h=$Hilb$^2(\mathbb{P}^1\times\mathbb{P}^1)$ and
determine recursively most
of the Big Quantum Cohomology. We prove a relationship between 
the invariants 
so obtained and the enumerative geometry of hyperelliptic curves in
 $\mathbb{P}^1\times\mathbb{P}^1$. This extends the results obtained by
Graber \cite{gr} for Hilb$^2(\mathbb{P}^2)$ and hyperelliptic curves in 
$\mathbb{P}^2$.}
\setcounter{section}{-1}

\section{Introduction}

The Gromov Witten invariants for a smooth projective complex variety $X$ count 
the (virtual) number of curves on $X$ satisfying some incident
conditions. Those related to the genus zero (virtual) curves are equivalent to
the Big Quantum Cohomology $QH^*$ of $X$; 
in particular the $3$-point invariants give the Small Quantum Cohomology
$QH^*_s$. If $X$ is convex, i.e. $H^1(\mathbb{P}^1,f^*(T_X))=0$ for all
genus zero stable maps $f:\mathbb{P}^1\rightarrow X$, the invariants are
enumerative since the moduli space is smooth of the 
expected dimension. In most cases where the invariants have been computed,
$QH^*_s$ is explicitely given by generators and relations; other genus
zero invariants are recursively determined.\\
Among the non-convex varieties whose $QH^*$ is known there is the Hilbert
scheme Hilb$^2(\mathbb{P}^2)$, studied by Graber in \cite{gr}. There he
computes $QH^*_s$ explicitely and $QH^*$ recursively using the First
Reconstruction Theorem FRT \cite{km}. He then relates the invariants to the
enumerative geometry of hyperelliptic curves in $\mathbb{P}^2$.\\
In this paper, we study the analogous problem for
$Q=\mathbb{P}^1\times\mathbb{P}^1$. The results we obtain are similar. The
main differences are: (1) $H^*($Hilb$^2(Q))$ is not generated by the divisor
classes hence a straightforward application of FRT is not enough to determine
all the invariants; (2) the group of automorphisms $Aut(Q)$ has four orbits on
Hilb$^2(Q)$, while $Aut(\mathbb{P}^2)$ has only two on Hilb$^2(\mathbb{P}^2)$;
(3) we have to be careful about intersection properties of curves, in
particular the Position Lemma \ref{lem:pos} does not apply directly. \\
In the first section we collect some basic facts about $\h=$Hilb$^2(Q)$.\\
Section two contains a detailed study of the deformation theory of certain
stable maps to $\h$, which is then used to compute some virtual classes and
therefore some invariants.\\
In section three the Small Quantum Cohomology of $\h$ is explicitely computed
and an algorithm is given; it determines most of the other genus zero GW
invariants. The methods used are a combination of classical enumerative
geometry, application of WDVV and the results of \S\ref{sec:deformazione}.\\
Section four contains the enumerative applications, i.e. an explicit relation
(Theorem \ref{hyp2}) between genus zero GW invariants and enumerative geometry 
of hyperelliptic curves on $Q$, which is completely analogous to that proven
for $\mathbb{P}^2$ by Graber. The proof is however considerably more
complicated.\\

\noindent
\textbf{Acknowledgements:} I would like to thank Professor Barbara Fantechi
for having introduced me to the topic of this paper and for helpful
discussions about it. \\
I am grateful to Professor Angelo Vistoli because during my visit to the
Dipartimento di Matematica at the Universit\`a di Bologna I could learn
a lot from him about the fascinating world of stacks.\\
A special thank goes to Professor Andr\'e Hirschowitz who made possible my 
visit to the Laboratoire J.A. Dieudonn\'e in Nice and to Joachim Kock for 
helpful discussions at the beginning of my work.\\
Finally I thank a lot SISSA and ICTP for the stimulating environment these
institutions provided me for the last months of my Ph.D.
\section{The Hilbert scheme Hilb$^2(\mathbb{P}^1\times\mathbb{P}^1)$}

In this section we fix notations and collect some results  
on the Hilbert scheme
 $\h:=$Hilb$^2(\mathbb{P}^1\times\mathbb{P}^1)$, 
whose points parametrize 0-dimensional length-2  
closed subschemes $Z$ of $\mathbb{P}^1\times\mathbb{P}^1$.\\

\noindent
\textbf{Notations and conventions:} we work over $\mathbb{C}$ and 
we identify
the variety $\mathbb{P}^1\times\mathbb{P}^1$ with its image under the Segre
embedding, i.e. the
smooth quadric $Q$ in $\mathbb{P}^3$. We have two rulings on $Q$, if
$q_1,q_2$ are the two projections on $\mathbb{P}^1$, the fibers of
 $q_1$ form the first ruling and those of $q_2$ the second one.\\
We consider Chow rings with $\mathbb{Q}$-coefficients. All the varieties in 
this section have a cellular decomposition, hence 
we write $A^*(X)$ for the Chow group of $X$ with $\mathbb{Q}$-coefficients. By 
\cite{ful} Example 19.1.11, $A^*(X)\cong H^{2*}(X)$. In particular we can identify them.\\
Let $\e$ be the sheaf of sections of a vector bundle $E$ ,
we denote by $\mathbb{P}(E)$ the projective 
bundle Proj(Sym$\e$). 
Geometrically, points of $\p{E}$ correspond to hyperplanes in the fibers
of $E$.\\
We indicate a non-reduced $0$-dimensional subscheme $Z$ of length $2$ of $Q$ 
as a pair $(p,v)$ where $p\in Q$ is the support of $Z$ and $v\in\p{T_{Q,p}}$ 
is a direction. We call it a non-reduced point of $\h$.
\subsection{Two geometrical descriptions of Hilb$^2(\mathbb{P}^1\times\mathbb{P}^1)$}\label{sec:descriptions}

There are two possible geometric descriptions 
of $\h$. The first is the standard one as 
a desingularization of the second symmetric product 
Sym$^2(\mathbb{P}^1\times\mathbb{P}^1)$ 
(see \cite{Fo}). Let $U$ be the product $Q\times Q$, $pr_1,pr_2$ the 
two projections, 
$\tilde{U}$ the blowup of $U$ along the diagonal $\delta\subseteq U$.
The group $\mathbb{Z}_2$ acts on $U$ fixing
$\delta$, so there is an induced action on the blowup $\tilde{U}$. The
 Hilbert scheme $\h$ is the quotient scheme $\tilde{U}/\mathbb{Z}_2$, hence it
 is smooth, projective, irreducible and $4$-dimensional. \\ 
We have the following diagram:
\begin{diagram}[height=0.5cm,width=1.0cm]
\tilde{\delta}  &  \rTo^j  & \tilde{U} & \rTo^{\theta} & \h\\
\dTo^{bl|_{\delta}} & & \dTo_{bl} & & \\
\delta &\rTo_i  & U & \pile{\rTo^{pr_1} \\ \rTo_{pr_2}} & Q \\
\end{diagram}
\noindent
with $i,j$ the natural inclusions, $bl$ the blowup map and  
$\theta$ the quotient map. It induces an isomorphism of 
$\mathbb{Q}$-algebras 
$\theta^*:\ah{*}\rightarrow A^*(\tu)^{\mathbb{Z}_2}$ which does not respect 
the degree. We denote by
$\D\cong\p{T_Q}$ the image in $\h$ of the exceptional divisor 
$\tilde{\delta}$.\\
Let $h_1, h_2$ be the cycle classes of the two rulings on $Q$. Then
$h_0=[Q]$, $h_1$, $h_2$, $h_3:=h_1h_2$ is a basis of 
$\aq{*}$ and $h_r\otimes h_s$, with $0\leq r,s \leq 3$, a basis of
 $A^*(U)$. Let $\xi$ be the class of the exceptional divisor $\tilde{\delta}$
 in $\tu$. Hence as 
a $\mathbb{Q}$-algebra $\ah{*}$ is generated by 
$~\xi,~T_1:=h_1\otimes 1+1\otimes h_1,~T_2:=h_2\otimes 1+1\otimes h_2,~
T_4:=h_3\otimes 1+1\otimes h_3$.

\vspace{0.5cm}
\noindent
The Hilbert scheme $\h$ can also be viewed as a blow up of the smooth 
projective 4-dimensional Grassmannian 
$\g:=$~Grass$(2,4)$ of lines in $\mathbb{P}^3$. In fact
there are two special lines  $\w{1},\w{2}\subseteq\g$ which are 
disjoint. A point $l_i\in\w{i}$ represents a line on the $i$-th ruling of $Q$, 
$i=1,2$. Denote by $W$ the disjoint union of these special lines,
i.e. $W=\{l\in\g:l\subseteq Q\}$. There exists a surjective
morphism $\varphi:\h\rightarrow \g$ defined by mapping a point $Z\in\h$ to its
associated line $l_Z$. 
\begin{theorem}\label{thm:descrizione2}
The Hilbert scheme $\h$ is isomorphic to the blow up of 
the Grassmannian $\g$ along $W$.
\end{theorem}
\begin{proof}
The morphism $\varphi$ is birational, its inverse is a morphism 
defined by $\varphi^{-1}(r)=r\cap Q$, for all $r\in\g-W$.\\
If $r\in W$ then the inverse image $\varphi^{-1}(r)$ is
Sym$^2(r)\cong\mathbb{P}^2$, so that $\varphi^{-1}(W)$ is a Cartier divisor in 
$\h$. Hence we have a commutative diagram:
\begin{displaymath}
\xymatrix{
\h\ar[r]^{\exists !~\alpha\;\;\;}\ar[rd]_{\varphi}& 
                                                     Bl_W \g \ar[d]^{\rho}\\
                                               & \g  }
\end{displaymath}
\noindent
where $\rho$ is the blowup morphism. Since both $\h$ and $Bl_W\g$ are smooth, 
$\alpha$ is an isomorphism if and only if it is bijective. It is obviously
bijective on $\g-W$. By explicit calculations it can be verified
that on the exceptional locus the generic fiber of $\alpha$ is a point, i.e.
$\alpha$ is a bijection. 
\end{proof}
\noindent
Let $\sig{1,0}\in\ag{*}$ be the Schubert cycle of
points $l\in\g$ intersecting a given line $r\subseteq\mathbb{P}^3$. We set 
$T_3:=\varphi^*(\sig{1,0})$.
\subsection{The cone of effective curves}\label{sec:effcone}

We define three $3$-codimensional cycle classes and show that they generate 
the cone of effective curves in $\h$.\\

\noindent
Fix a point $l_1\in\w{1}$ and let $C(l_1)$ be a line in the plane Sym$^2(l_1)$.
Note that all the points $Z$ of $\h$ 
contained in $C(l_1)$ are such that $\supp{Z}\subseteq l_1$. We denote
by $\cuno$ the corresponding cycle class in $\ah{3}$. 
We define the class $\cdue$ analogously.
Fix a point $p_0\in Q$ and consider the line 
$C(p_0)=\p{T_{Q,p_0}}$. 
Let $\ctre$ be the corresponding cycle class in $\ah{3}$. 
Note that for all $Z\in C(p_0)$ we have 
$\supp{Z}=p_0$. 
The curves $C(l_1),~C(l_2),~C(p_0)$ are effective in $\h$.
\begin{proposition}\label{prop:cono-eff}
An effective curve in $\h$ is of class $a\cuno+b\cdue+c\ctre$ with 
$a,b,c\geq 0$. 
\end{proposition} 
\begin{proof}
The linear systems associated to $T_1,T_2,T_3$ (see \S 
\ref{sec:descriptions}) are base-points-free. Since: 
$$C_1\cdot T_2=1, \;C_2\cdot T_1=1,\;F\cdot T_3=1$$ 
\noindent
and all other 
possible intersections give zero, an effective curve
in $\h$ is of class $aC_1+bC_2+c\ctre$ with $a,b,c\geq 0$.
\end{proof}
\noindent
We will write $(a,b,c)$ for the class $a\cuno+b\cdue+c\ctre$.
\begin{rem}\label{rem:c-h}
$\h$ is the blowup of $\g$ along $W$,
so by \cite{gh} p.608 we have $c_1(T_{\h})=2(T_1+T_2)$.
\end{rem}
\begin{corollary}\label{expdim}
The expected dimension $\ed{\h}$ of $\mbar{0,n}(\h,(a,b,c))$ 
is given by the formula: 
\begin{displaymath} 
\textrm{exp.dim}~\mbar{0,n}(\h,(a,b,c))=\ed{\h}=2a+2b+1+n
\end{displaymath}
\end{corollary}
\begin{proof}
It follows from \ref{rem:c-h} and the general formula for the expected 
dimension of a moduli space of stable maps (see \cite{beh}).
\end{proof}
\subsection{A good $\mathbb{Q}$-basis for $\ah{*}$}\label{sec:good-basis}

By sections \ref{sec:descriptions} and \ref{sec:effcone} 
we know that $\ah{*}$ can be generated by $T_1$, $T_2$, $T_3$ and $T_4$. 
Note that
$T_4$  can be represented by the cycle class: 
$$\G{p}=\{[Z]\in\h:p\in\supp{Z},~p\in Q\; \textrm{given 
point}\}$$ 
\noindent
the inverse image of the cycle $\sig{2,0}(p)$ in $\g$ via the
blowup map $\varphi$.

\vspace{5pt}
\noindent
We complete once for all $T_1,T_2,T_3, T_4$ to a  
$\mathbb{Q}$-basis for $\ah{*}$ by adding: $T_0=[\h]$, 
~$T_5=T_1T_2$, ~$T_6=T_1^2$, ~$T_7=T_2^2$, ~$T_8=T_1T_3$, ~$T_9=T_2T_3$,~
$T_{10}=\cdue+\ctre$, ~$T_{11}=\cuno+\ctre$,
~$T_{12}=\cuno+\cdue+\ctre$ and $T_{13}$ the class of a point.
\begin{rem}\label{rem:sym}
Let $\iota$ be the involution of $Q=\mathbb{P}^1\times\mathbb{P}^1$ defined by 
$\iota(p,q)=(q,p)$. The induced involution on $\h$, also denoted by $\iota$, 
interchanges $T_1$ and $T_2$, $T_6$ and $T_7$, 
$T_8$ and $T_9$, $T_{10}$ and $T_{11}$, and leaves the other $T_i$'s invariant.
\end{rem}
\noindent
\subsection{The action of $Aut(Q)$ on $\h$}\label{sec:aut}

Let $\mathcal{A}_0=PGL(2)\times PGL(2)$ be the connected component containing 
the identity in $\mathcal{A}=Aut(Q)$. Then 
$\mathcal{A}=\mathcal{A}_0\sqcup\iota\mathcal{A}_0$.
There are four orbits on $\h$ with  respect to the
$\mathcal{A}$-action and we can give a description of all of them:
\begin{displaymath}
\begin{array}{l}
\mo_4=\{Z\in\h: \supp{Z}=\{p,q\},p\neq q, l_Z\nsubseteq Q\}\\
\mo_3=\{Z\in\h: \supp{Z}=\{p,q\},p\neq q, l_Z\subseteq Q\}\\
\D_3=\{Z\in\h: \supp{Z}=p, l_Z\nsubseteq Q \}\\
\D_2=\{Z\in\h: \supp{Z}=p, l_Z\subseteq Q \}
\end{array}
\end{displaymath}
\noindent
Here indexes are chosen equal to the dimensions of the orbits.\\
The closed orbit $\D_2$ is the disjoint union of two closed subvarieties
$\D_2^i$, $i=1,2$, where $Z\in\D_2^i$ if $l_Z\in W_i$.\\ 
Note that $\overline{\D}_3=\D$ is the divisor of
non-reduced points, i.e. the image of $\tilde{\delta}$.\\
The orbit $\mo_3$ is the disjoint union $\mo_3^1\sqcup\mo_3^2$ where
$$\mo_3^i=\{Z\in\h: \supp{Z}=\{p,q\},p\neq q,~l_Z\in\w{i}\}$$ 
\noindent
In particular  
the closures  $\omo_3^1$, $\omo_3^2$ are the two exceptional divisors $\wt{1}$,
$\wt{2}$ respectively, of
the blowup map $\varphi:\h\rightarrow\g$.\\ 
Finally the orbit $\Sigma_4$ is open and dense in $\h$.\\
$\h$ is an \emph{almost-homogeneuos space} since it has a
finite number of orbits for the $\mathcal{A}$-action, hance we can use Graber's
 Position Lemma. 
\begin{lemma}\label{lem:pos}
\emph{\textbf{(Position Lemma-\cite{gr} Lem.2.5)}}
Let $A$ be a smooth, almost-homogeneous space under 
the action of an integral group $G$, $f:B\rightarrow A$ a morphism with 
$B$ smooth.
Let $\Gamma$ be a smooth cycle on $A$ which intersects the stratification 
properly, and $\Gamma_{\textrm{reg}}$ be the locus in $\Gamma$ 
where the intersection with the stratification is transversal. Then:
\begin{enumerate}
\item for a generic $g\in G$, $f^{-1}(g\Gamma)$ is of pure dimension equal 
to the expected one;
\item the open set (possibly empty) $f^{-1}(g\Gamma_{\textrm{reg}})$ 
is smooth. 
\end{enumerate}
\end{lemma}
\begin{rem}\label{rem:weak-pos}
If in the hypotheses of \ref{lem:pos} we do not assume $B$ smooth but only pure
dimensional we can consider its desingularization $\nu:\tilde{B}\rightarrow
B$. Then by applying the Position Lemma to the composition map 
$\tilde{f}:\tilde{B}\rightarrow A$ we get that
$\cod(\tilde{f}^{-1}(g\Gamma)\subseteq\tilde{B})$ is the expected 
one, i.e. equal to $\cod(g\Gamma\subseteq A)$. Since:
$$\cod(\tilde{f}^{-1}(g\Gamma)\subseteq\tilde{B})\leq
\cod(f^{-1}(g\Gamma)\subseteq B)$$ 
\noindent
we have that \ref{lem:pos}-1) holds with the inequality 
$$\cod(f^{-1}(g\Gamma)\subseteq B)\geq\cod(g\Gamma\subseteq A)$$ 
\end{rem}
\begin{rem}\label{rem:pos-0}
Note that $\mathcal{A}$ is not integral, so we will apply the Position
Lemma \ref{lem:pos} to $G=\mathcal{A}_0$. 
\end{rem}
\begin{rem}\label{rem:inters}
For any $p\in Q$ the cycle $\G{p}$ intersects the stratification properly. In
fact $\G{p}\cap\mo_4\cong Q-(l_1(p)\cup l_2(p))$ is obviously a proper
intersection and $\G{p}\cap \D_2=\{(p,T_{l_1(p),p}),~(p, T_{l_2(p),p})\}$ is
$0$-dimensional. Since these intersections are non-empty, it is also satisfied
$\G{p}\nsubseteq\mo_3\sqcup\D_3$. 
\end{rem}
\noindent
We set $\G{p}_{\textrm{reg}}$ to be the locus of $\G{p}$ where the
intersection with the stratification is transversal. 
\begin{lemma}\label{lem:trans}
Given a point $p\in Q$, $\G{p}_{\textrm{reg}}$ is the open subset 
of $\Gamma(p)$ of points with reduced support.
\end{lemma}
\begin{proof}
We first prove that $\D_2^k\cap\G{p}$, $k=1,2$, is not transversal. $\D_2^k$ is
 the pullback of the diagonal $\D$ via the inclusion map 
$j_k:\wt{k} \hookrightarrow\h$, hence it is a divisor in $\wt{k}$. By the 
projection formula we obtain $\D_2^k=4l_k-2\xi_k$. It is easy to verify that
 $\G{p}$ intersects $\D_2^k$ only in one point, but $T_4\cdot (j_k)_*\D_2^k=2$,
 this means the intersection is not transversal.\\ 
Now consider the closed immersion $f:Q\rightarrow Q\times Q$, defined by 
$f(q)=(p,q)$.
Let $\theta$ be the quotient map defined in \S \ref{sec:descriptions}.
There exists a unique induced closed immersion 
$\tilde{f}:Bl_p Q=\theta^{-1}(\G{p})\to\tu$, by \cite{har}, Chap.II Cor.7.15.
Choosing local coordinates on $\h$ and $\tu$ it is easy to see that for
each $(p,v)\in\theta^{-1}(\D)$ the image $d\theta_p(T_{(p,v)}\tu)$ is 
contained into $T_{\theta (p,v)}\D$. As $\tilde{f}$ is a closed immersion 
and $\D_3$ is open dense in $\D$, it follows that $\G{p}$ does not intersect 
$\D_3$ transversally.\\
In order to study 
the differential of the map $\G{p}\rightarrow\mo_3$ it is enough 
to restrict it to the divisor $\wt{k}$ and to study the 
differential of $Q-\{p\}\rightarrow Q\times Q-\delta$,
$q\mapsto (p,q)$. As $Q=\mathbb{P}^1\times\mathbb{P}^1$ we can choose affine 
coordinates on both $\mathbb{P}^1$'s so that $p=(p_1,p_2)$ and the above map 
becomes:
\begin{displaymath}
\begin{array}{rcl}
\mathbb{A}^2-\{(p_1,p_2)\} & \rTo & \mathbb{A}^4\\
(q_1,q_2)& \mapsto &(p_1,p_2,q_1,q_2)
\end{array}
\end{displaymath}
Denoting by $x_1,x_2,y_1,y_2$ the coordinates on $\mathbb{A}^4$, the tangent 
space $T_{(p,q)}\wt{k}$ is the $3$-dimensional 
affine space defined by the equation $x_k-y_k=0$. In these coordinates
 $\G{p}-\D$ is the set
 $\{(x_1,x_2,y_1,y_2): x_1=p_1,x_2=p_2\}$ so its tangent space at $(p,q)$
 is the $2$-dimensional affine 
space defined by the equations $x_1=0, x_2=0$. Hence for each 
$(p,q)\in\G{p}\cap\wt{k}-\D$, the space $T_{(p,q)}(\G{p}-\D)$ is not contained 
in $T_{(p,q)}\wt{k}$, that is to say $\G{p}$ intersects $\mo_3$ transversally.
\end{proof}
\subsection{The locus $\D$ of non-reduced points of $\h$}\label{sec:delta}

We will refer to the locus $\D$ of non-reduced points of $\h$ as the 
\emph{diagonal} of $\h$. Its class in $\ah{1}$ is $2(T_1+T_2-T_3)$. Given a 
line $l$ in $Q$, $\D$  intersects the fiber $\varphi^{-1}(l)\cong~
$Sym$^2(l)$ in a smooth conic.\\    
The natural map $s:\D\rightarrow Q$
 is defined by mapping a non reduced point to its support 
so we will call it the \emph{support map}. 
\begin{proposition}\label{prop:pic-delta}
Let $i:\D\rightarrow\h$ be the inclusion. Then:
\begin{gather*}
\textrm{Pic}(\D)=\langle \frac{1}{2}T_1,\frac{1}{2}T_2,T_3 \rangle
\end{gather*}
\end{proposition}
\begin{proof}
Let $T_j=i^*T_j$, by abuse of notation.
The thesis follows from the equality $Pic(\p{T_Q})=\langle s^*Pic(Q), 
\mathcal{O}(1)\rangle$.
\end{proof}
\begin{proposition}\label{prop:eff-delta}
The effective curves in
$\h$ which are contained into $\D$ are of class $(a,b,c)$ with $a,b,c\geq 0$
and $a,b$ even.
\end{proposition}
\begin{proof}
Let $C\subseteq\D$ be an effective curve of class $(\alpha,\beta,\gamma)$, 
then $i_*C$ is an effective
curve in $\h$ of class $(a,b,c)$ for some non negative integers $a,b,c$.
By the projection formula, deg$_{\D}\frac{1}{2}T_1\cdot C=\frac{a}{2}$ 
is an integer number equal to $\alpha$, hence $a$ is even. 
The same is true for $b$, by symmetry.
\end{proof}
\begin{rem}\label{rem:c-d} 
By the adjunction formula and \ref{rem:c-h} we get 
$c_1(T_{\D})=2T_3$.
\end{rem}
\subsection{The divisor $\mo$}\label{sec:sigma}

Let $\mo$ be given by the disjoint union $\wt{1}\sqcup\wt{2}$ of the two 
exceptional divisors of the blowup map $\varphi$;
as an element of $\ah{1}$ it is the class $2T_3-T_1-T_2$.\\
Note that $\wt{1}$ is isomorphic to $\mathbb{P}^2\times\mathbb{P}^1$  
because it is the
relative Hilbert scheme $Hilb^2(Q/\mathbb{P}^1)$. 
Let $\pi_i:\wt{1}\rightarrow \mathbb{P}^i$ be the natural projection. Then 
Pic$(\wt{1})$ is generated by $L_1=\pi^*_1\mathcal{O}(1)$, 
$L_2=\pi^*_2\mathcal{O}(1)$ and $A_1^{eff}(\wt{1})$ is generated by 
$A_1=[\mathbb{P}^1\times pt]$ and $A_2=[pt\times\mathbb{P}^1]$. 
\begin{lemma}\label{lem:j1}
Let $j_1:\wt{1}\rightarrow\h$ be the inclusion. Then 
$j_1^*T_1=2L_2$, $j_1^*T_2=L_1$, $j_1^*T_3=2L_2$.
\end{lemma}
\begin{proof}
It is enough to compute deg$~T_i|_{A_j}$. 
Fix $l_2^{'}, l_2^{''}\in W_{2}$ and consider the curve $C=\{Z:
\supp{Z}=\{p,q\},~\exists~l_1\in\w{1}~\textrm{with}~p=l_1\cap
l_2^{'},~q=l_1\cap l_2^{''}\}$. Then $A_2=[C]$. 
It is easy to see that $j_{1*}(A_1)=\cuno$. Hence it is elementary to verify 
the claim.
\end{proof}
\begin{corollary}\label{rem:c-s}
$c_1(T_{\mo})=3(T_1+T_2)-2T_3$.
\end{corollary}
\begin{proposition}\label{prop:eff-w1}
An effective curve in $\h$ which is contained into $\wt{1}$ is of class
$(a,b,c)$ with $a,b,c\geq 0$ and $b=c$ even. 
\end{proposition}
\begin{proof} We know that $j_{1*}(A_1)=\cuno$, hence it is enough to prove
  that $j_{1*}(A_2)=2\cdue+2\ctre$. This follows by Lemma \ref{lem:j1} and the
  projection formula. 
\end{proof}
\subsection{Description of some effective curves}\label{sec:curves}

We describe all the effective connected curves in some cycle classes in 
$A_1(\h)$. For more details see \cite{p} \S 1.7.
In the following sections we will 
make explicit calculations on the moduli spaces of stable maps involving such
curves. \\
 
\noindent
$\mathbf{Notations:}$ If $p\in Q$ is a point we will denote by $l_i(p)$ 
the unique line of the
$i$-ruling on $Q$ going through $p$.\\
We will use Propositions \ref{prop:eff-delta} and \ref{prop:eff-w1} without
explicit reference throughout.
\subsection*{Curves of class $(0,0,c)$}

A curve of class $(0,0,c)$ is contained in $\D$: it is $c$ times a fiber 
$C(p)$ of the support map over some $p\in Q$.
\subsection*{Curves of class $(1,0,c)$, $(0,1,c)$}

Since the classes $(1,0,c),(0,1,c)$ are symmetric under the involution
 we can analyse only one of them. We choose $(1,0,c)$.\\
Let $\tilde{\varphi}:\h\rightarrow Bl_{W_1}\g$ be the natural map.
A curve of class $\beta=(1,0,0)$ is a line in a fiber
 $\tilde{\varphi}^{-1}(l_1)=$Sym$^2(l_1)$, for some $l_1\in W_1$.
We will denote it by $C(l_1)$, (see \S \ref{sec:effcone}).\\
An irreducible curve of class $(1,0,1)$ has the form 
$C(p_1,l_1)=\{Z\in\h: \supp{Z}=(p_1,q), q\in l_1\}$, where $p_1\in Q$ is a
 fixed point and $l_1\in W_1$ a fixed line such that $p_1\notin l_1$.  
A reducible curve of this class is the union of two irreducible effective 
components  $C(l_1)\cup C(p)$ for some $l_1\in W_1$ and $p\in l_1$, with $p\in
C(l_1)\cap \D|_{\varphi^{-1}(l_1)}$. Note that it is contained into
 $\D\cup\mo$.\\
For $c\geq 2$ there are only reducible curves of class $(1,0,c)$: they are
 entirely contained into $\D\cup\mo$ with support 
$C(l_1)\cup C(p)\cup C(q)$ or 
$C(l_1)\cup C(p)$ for some $l_1\in W_1$ and  $p,q$ points in
$C(l_1)\cap\D|_{\varphi^{-1}(l_1)}$.
\subsection*{Curves of class $(1,1,c)$, $c\leq 1$ }

Connected curves of class $(1,1,0)$ do not exist.\\
Let $C$ be a reducible curve of type $(1,1,1)$. We have three possible
  decompositions:\\
- $C(l_1(p))\cup C(l_2(p))\cup C(p)$ with $p$ a point of $Q$;\\
- $C(p,l_1)\cup C(l_2(p))$ for a given line $l_1$ and a given point $p\in Q$,
with $C(l_2(p))$ a line in
  Hilb$^2(l_2(p))$ passing through $(p,q),~q=l_1\cap l_2(p)$;\\
- $C(p,l_2)\cup C(l_1(p))$ symmetrically. \\
We have two possible families of irreducible curves of 
class $(1,1,1)$. 
Fix a plane $\Lambda\subseteq\mathbb{P}^3$ and a generic point 
$q\in \Lambda$, $q\notin Q$. 
If $\Lambda$ is generic a curve of such a class is 
a line $\Ll$ in Hilb$^2(\Lambda\cap Q)$, whose points are 
the closed subschemes $Z$ such that $\supp Z\subseteq (\Lambda\cap Q)$, 
$q\in l_Z$. Otherwise the irreducible curve is determined by choosing
$\Lambda$ 
tangent to $Q$ at a point $p$ and $q\in \Lambda$ such that 
$q\notin \Lambda\cap Q$. 
Its points are the closed subschemes $Z$ such that $\supp Z\cap
l_1(p)\neq\emptyset$, $\supp Z\cap l_2(p)\neq\emptyset$ and $q\in l_Z$. Note
that such a curve has only $4$ moduli, while the expected dimension is $5$.
\begin{rem}\label{rem:gen-curve}
Note that irreducible curves $C(p_1,l_1)$, $C(p_2,l_2)$,
$\Lambda(l)$, of class $T_{11}$, $T_{10}$ and  
$T_{12}$ respectively, intersect the stratification properly.
\end{rem}
\section{Virtual Fundamental Classes}\label{sec:deformazione}

This section presents some
 results about the way of computing some GW invariants we will 
need in the following. In particular we calculate the virtual fundamental
 class of two moduli spaces of stable maps on $\h$.   
\subsection{Deformation theory on $\mbar{0,n}(\h,\beta)$}\label{sec:def.th.}

We recall the following fundamental result (see \cite{k} Thm.II.1.7):
\begin{theorem}\label{thm:smooth}
If $\mu:C\rightarrow\h$ is a $n$-pointed stable map and 
$H^1(C,\mu^*T_\h)=0$, then the forgetful morphism
$\eta:\mbar{g,n}(\h,\beta)\rightarrow \mathfrak{M}_{g,n}$ is smooth at 
the point $[C,x_1,\ldots,x_n,\mu]$. 
\end{theorem}  
\begin{rem}\label{rem:convex}
A smooth variety $X$ is called \emph{convex} if $H^1(\mathbb{P}^1,f^*T_X)=0$ 
for all genus zero stable maps $f:\mathbb{P}^1\rightarrow X$. 
If $X$ is convex then $H^1(C,f^*T_X)=0$ for all maps $f:C\rightarrow X$, $C$ a
genus zero rational curve. Hence the moduli space 
$\mbar{0,n}(X,\beta)$ is smooth of dimension equal to the expected one, 
\cite{al} I.3. Note that every homogeneous variety is convex.
\end{rem}
\noindent
In order to compute some GW invariants, we need to control the smoothness 
of the moduli space $\mbar{0,n}(\h,\beta)$. 
$\h$ is not convex but it is an almost-homogeneous space under the action of
$\mathcal{A}$.
\begin{theorem}\label{lem:sigma-convex}
a) $\mo$ is convex hence $\mbar{0,n}(\mo,(a,b,c))$ 
is smooth of the expected dimension $\ds{\mo}=3(a+b)-2c+n$.\\
b) $\D$ is convex hence $\mbar{0,n}(\D,(a,b,c))$ is smooth of
the expected dimension $\ds{\D}=2c+n$.
\end{theorem}
\begin{proof}
Recall that $\mo=\wt{1}\sqcup\wt{2}$.
$\wt{k}$ is isomorphic to $\mathbb{P}^2\times\mathbb{P}^1$ therefore it is
homogeneous, hence convex. The formula for 
$\ds{\mo}$ follows from \ref{rem:convex} and \ref{rem:c-s}.
This proves statement a).\\ 
The support map $s:\D\rightarrow Q$ gives the exact sequence:
\begin{displaymath}
0\rightarrow T_{\D/Q}\rightarrow T_{\D}\rightarrow s^*T_{Q}\rightarrow 0 
\end{displaymath}
\noindent
Let $\mu:\mathbb{P}^1\rightarrow \D$ be a stable map, then we have to prove 
that $H^1(\mathbb{P}^1,\mu^*T_{\D/Q})$ vanishes, since $Q$ is homogeneous.
The generators of the cone of effective 
curves in $\D$ are such that the degree of $T_{\D/Q}$ restricted to each of
them is non-negative, so deg$~\mu^*T_{\D/Q}\geq 0$ and
 $H^1(\mathbb{P}^1,\mu^*T_{\D/Q})=0$.
As before the formula for $\ds{\D}$ follows from 
\ref{rem:convex} and \ref{rem:c-d}. This concludes the proof.
\end{proof}
\begin{theorem}\label{thm:dim-h}
If $\mu:C\rightarrow \h$ is a stable map from a genus $0$ curve such that
 no component of $C$ is mapped entirely into $\D\cup\mo$, 
then the moduli space 
$\mbar{0,0}(\h,\beta)$ is smooth at $[C, \mu]$ of the expected dimension.
\end{theorem}
\begin{proof}
$\h-(\D\cup\mo)$ is $\mo_4$, 
the open dense orbit for the action on $\h$ induced by 
$\mathcal{A}$. The action on $\mo_4$ is
transitive, so  we can say that $T_{\h}$ is generically generated 
by global sections on $\h$. Let $\mu:C\rightarrow \h$ 
be as in the hypotesis, then $\mu^*T_{\h}$ 
is generically generated by global sections on $C$. This means that 
$H^1(C,\mu^*T_{\h})=0$ and  the moduli space $\mbar{0,0}(\h,\beta)$
is smooth at $[C, \mu]$ of the expected dimension by \ref{thm:smooth}. 
\end{proof} 
\subsection{The moduli space $\mbar{0,0}(\h,(0,0,c))$}\label{sec:c}

Here and in the following section we prove some results on the virtual 
fundamental class of two moduli spaces which we will use later on to 
make explicit calculations. \\

\noindent
The virtual fundamental class is defined on $\mbar{g,n}(X,\beta)$ by using 
a perfect obstruction theory (in the sense of \cite{bf}, \cite{beh}) or, 
equivalently, a tangent-obstruction complex (as in \cite{lt}).
\begin{proposition}\label{prop:obs-bdle}
On the open locus where the moduli stack is smooth, the virtual fundamental 
class can be computed by integrating the top Chern class of the 
\emph{obstruction bundle}.  
\end{proposition}
\begin{proof}
\cite{bf}, Proposition 5.6.
\end{proof}
\noindent
We will always denote the obstruction bundle by $\e$.\\
If $u:\mathcal{C}\rightarrow\mbar{g,n}(X,\beta)$ is the universal curve with 
universal map $f:\mathcal{C}\rightarrow X$, then  
$\e=\ext_u^2(f^*\Omega_X\rightarrow\Omega_u,\mathcal{O}_{\mathcal{C}})$, 
whose fiber over a point $[C,x_i,\mu]$ is $\ext^2(f^*\Omega_X\rightarrow
\Omega_C,\mathcal{O}_C)$. It follows that the obstruction bundle on 
$\mbar{g,n}(X,\beta)$ is the pullback of the obstruction bundle on 
$\mbar{g,0}(X,\beta)$.
\begin{proposition}\label{rel-obs}
If the moduli space $\mbar{g,n}(X,\beta)$ is smooth over the Artin 
stack $\mathfrak{M}_{g,n}$ then $\e=\R^1\pi_* ev^* T_X$ is a relative 
perfect obstraction theory for $\mbar{g,0}(X,\beta)$,   
with $ev$ the 
usual evaluation map and $\pi$ the flat morphism forgetting the $n$ marked 
points and stabilizing.
\end{proposition}
\begin{proof}
\cite{beh}, Proposition 5.
\end{proof}
\noindent
Let us consider the moduli space $\mbar{0,0}(\h,(0,0,c))$.\\ 
For $c\geq 1$, a curve of class $(0,0,c)$ in $\h$ is represented by a 
$c$-sheeted cover of $\mathbb{P}^1$ and it is contained into $\D$ which is
convex. 
Then the moduli space 
$\mbar{0,0}(\h,(0,0,c))$ is smooth of dimension $2c$ bigger than the expected
one, $\ed{\h}=1$. The obstruction bundle $\e=\R^1\pi_*(ev^*T_\h)$ has rank 
 $2c-1$ and its stalk at the point $[C,\mu]$ is 
$H^1(C,\mu^*T_\h)= H^1(C,\mu^*\mathcal{O}_{\mathbb{P}^1}(-2))$. 
By \ref{rel-obs} the virtual fundamental class is given by the product:
\begin{displaymath}
[\mbar{0,0}(\h,(0,0,c))]^{vir}=[\mbar{0,0}(\h,(0,0,c))]\cdot c_{2c-1}(\e)
\end{displaymath}
\begin{proposition}\label{prop:ec}
Let $g:\mbar{0,0}(\h,(0,0,c))\rightarrow Q$ be the map defined by
$g([C, \mu])=\supp{\mu(C)}$. It holds: 
\begin{equation}\label{eq:ctop-gen}
c_{2c-1}(\e)=-g^*K_Q\cdot c_{2c-2}(\tilde{\e})
\end{equation}
\noindent 
where $\tilde{\e}$ is such that:
\begin{equation}\label{eq:c-cubo}
c_{2c-2}(\tilde{\e}|_{g^{-1}(p)})=\frac{1}{c^3}
\end{equation}
\noindent
for any point $p\in Q$.
\end{proposition}
\begin{proof} 
Let $\tilde{ev}:\mbar{0,1}(\h,(0,0,c))\rightarrow\D$
be the evaluation map into $\D$ such that the composition with the inclusion 
$\D\hookrightarrow\h$ gives the usual evaluation 
$ev:\mbar{0,1}(\h,(0,0,c))\rightarrow \h$. 
By \cite{lq} Lemma 3.2, 
$\e$ sits in the exact sequence:  
\begin{displaymath}
0\rightarrow g^*\mathcal{O}_Q(-K_Q)\rightarrow\e\rightarrow
\R^1\pi_*\tilde{ev}^*(s^*T_Q\otimes\mathcal{O}_\D(-1))=\tilde{\e}\rightarrow 0
\end{displaymath}
\noindent
Hence we get:
\begin{displaymath}
c_{2c-1}(\e)=-g^*K_Q\cdot c_{2c-2}(\tilde{\e})
\end{displaymath}
\noindent
Note that the inverse image $g^{-1}(p),~ p\in Q$, is isomorphic to 
$\mbar{0,0}(\mathbb{P}^1,c)$, with $\mathbb{P}^1\cong M_2(p)$ the
punctual Hilbert scheme of points on $Q$ at $p$.\\
 With respect to the diagram:
\begin{displaymath}
\begin{diagram}[height=0.8cm,width=1.0cm]
\mbar{0,1}(\mathbb{P}^1,c) & \rTo^{ev_1} & \mathbb{P}^1 \\
\dTo^{f} & & \\
\mbar{0,0}(\mathbb{P}^1,c) & &
\end{diagram}
\end{displaymath}
\noindent
the restriction $\tilde{\e}|_{g^{-1}(p)}$ is isomorphic to
$\R^1f_*ev_1^*(\mathcal{O}_{\mathbb{P}^1}(-1)\oplus
\mathcal{O}_{\mathbb{P}^1}(-1))$ \cite{lq} Rmk.3.1. By Theorem
3.2 in \cite{man}:
$$c_{2c-2}(\R^1f_*ev_1^*(\mathcal{O}_{\mathbb{P}^1}(-1)\oplus
\mathcal{O}_{\mathbb{P}^1}(-1)))=\frac{1}{c^3}$$
\noindent
This concludes the proof.
\end{proof}
\subsection{The moduli space $\mbar{0,0}(\h,(1,0,c))$}
\label{sec:m(c)}

In this section we want to study the moduli spaces $\mbar{0,0}(\h,(1,0,c))$ and
 compute (part of) their virtual fundamental classes. By \ref{expdim} they 
have expected dimension equal to $3$.\\

\noindent
The moduli space $\mbar{0,0}(\h,(1,0,0))$ is smooth of the expected dimension, 
because 
$(1,0,0)$ is the class of a curve contained into $\mo$ which is convex and
$\mbar{0,0}(\mo,(1,0,0))$ has the same expected dimension (by
\ref{lem:sigma-convex}).\\
We recall from section \ref{sec:curves} that the only irreducible curves of 
class $(1,0,c)$ must have $c=0$ or $1$; they are all smooth and rational, and 
if $c=1$ they are disjoint from $\D$. In $\mbar{0,0}(\h,(1,0,1))$ the open 
substack of maps having irreducible domain is also closed. 
It is easy to see (\cite{p} \S 2.6) that 
$\mbar{0,0}(\h,(1,0,1))$ is smooth of the expected dimension.\\
Irreducible curves of class $(0,0,c)$ 
necessarily have $c=1$, they are smooth, rational and contained in $\D$. Hence 
the domain of every stable map in $\mbar{0,0}(\h,(1,0,c))$ is reducible for 
$c\ge 2$. If $\mu: C\to\h$ is a stable map of class $(1,0,c)$ with $c\ge 2$,
$C$ has a unique component $C_0$ mapping isomorphically to a curve
of class $(1,0,0)$. This defines an $\mathcal{A}_0$-equivariant morphism:
\begin{displaymath}
\tau:\mbar{0,0}(\h,(1,0,c))\rTo\mbar{0,0}(\h,(1,0,0))
\end{displaymath}
\noindent
For $c\ge 1$ let \emph{$M(c)$} be the smooth closed substack $ev^{-1}(0)$ in
$\mbar{0,1}(\mathbb{P}^1,c)$. Let \emph{$M(0)$} be a point. Note that 
dim $M(c)=2c-2$ for $c\ge 1$ and that $M(1)$ is also a point.
\begin{lemma}\label{lem:descrizione}
The general fiber of $\tau$, i.e. over a curve intersecting $\D$ transversally
in two points, is isomorphic to:
$$\coprod\limits_{\stackrel{c_i\geq 0}{c_1+c_2=c\geq 2}}M(c_1)\times M(c_2)$$
\noindent
In particular $\tau$ is smooth over the open dense $\mathcal{A}_0$-orbit in
$\mbar{0,0}(\h,(1,0,0))$. 
\end{lemma}
\begin{proof}
Let $C$ be a general curve of class $(1,0,0)$. Since $C$ is fixed as well as 
its intersection points with the diagonal, 
the only moduli comes from the choice of the sheeted covers of the 
$(0,0,1)$-curves, i.e. curves in  
$M(c_i)$, $i=1,2$ with $c_1+c_2=c\geq 2$, such that the marked point mapping 
to the origin of the $(0,0,1)$-curve is in $C\cap\D$. 
\end{proof}
\noindent
\begin{rem}\label{rem:liscezza m(c)}
The composition of the inclusion: 
$$M(c_1)\times M(c_2)\rightarrow\mbar{0,0}(\h,(1,0,c))$$
\noindent
with the forgetful map 
$\mbar{0,0}(\h,(1,0,c))\rightarrow\mathfrak{M}_{0,0}$ is smooth
on its image which consists of the (smooth) locus of codimension $2$ 
parametrizing curves with two nodes $q_1, q_2$ if $c_1,c_2>0$, and the 
divisor parametrizing curves with a 
node in $q_1$ if $c_2=0$ or in $q_2$ if $c_1=0$.
\end{rem}
\begin{rem}\label{rem:dim-fibra}
 A general fiber of $\tau$ has expected dimension
equal to zero, since $\mbar{0,0}(\h,(1,0,c))$ has expected dimension 
$\ed{\h}=3$ and $\tau$ is smooth on the open dense orbit.\\
Its virtual fundamental class
$[\mbar{0,0}(\h,(1,0,c))]^{vir}\cdot\tau^*[C,f]$ is equal to 
the sum of the virtual fundamental
classes of all components. Moreover each of them must have expected dimension
equal to zero.  
\end{rem}
\noindent
To calculate the virtual fundamental class of a general fiber of $\tau$ 
we need to know the obstruction bundle $\e$ at one of its points.
 The following lemma gives a description of the space 
$H^1(D,\mu^*T_\h)$ which will permit us to express $\e$ as the cokernel of an
injection (see \ref{prop:e}). Fix $[D,\mu]$ a point in such a fiber:
\begin{displaymath}
\begin{array}{ll} 
\mu:D_0\rTo^{\cong} C(l_1) & \\
\mu:D_i\rTo^{c_i:1} C(p_i) & \\
\mu(q_i)=p_i\in C(l_1)\cap C(p_i) & 
\end{array}
\end{displaymath}
\noindent
We assume $c_i>0$, for $i=1,2$; the case with a $c_i$ equal to $0$ is 
similar but easier. 
\begin{lemma}\label{lem:l}
Let $\mathcal{L}_i$ be the invertible sheaf 
$\mu^*\mathcal{O}_{C(p_i)}(-2)$ of degree $-2c_i$, $i=1,2$. Then:
\begin{displaymath}
H^1(D,\mu^*T_\h)\cong H^1(D_1,\mathcal{L}_{1})\oplus H^1(D_2,\mathcal{L}_{2})  
\end{displaymath} 
\end{lemma}
\begin{proof}
We consider the exact sequence in cohomology:
\begin{displaymath}
\begin{array}{l}
H^0(T_{q_1}\oplus T_{q_2})\rightarrow H^1(D,\mu^*T_\h)
\rightarrow \bigoplus\limits_{i=1,2} H^1(D_i,\mu^*T_\h|_{D_i})\rightarrow 0
\end{array}
\end{displaymath} 
\noindent
The support map $s:\D\rightarrow Q$ is a $\mathbb{P}^1$-bundle, so the usual 
exact sequence:
\begin{displaymath}
0\rightarrow T_\D\rightarrow T_\h|_\D\rightarrow
\mathcal{N}_{\D/\h}\rightarrow 0
\end{displaymath}
\noindent
restricted to a fiber $l$  of $s$ gives
$\mathcal{N}_{\D/\h}|_l=\mathcal{O}_l(-2)$. Hence we get:
\begin{displaymath} 
H^1(D_i,\mu^*T_\h|_{D_i})=H^1(D_i,\mu^*\mathcal{O}_{C(p_i)}(-2))
\end{displaymath}
\noindent
and the above sequence becomes:
\begin{displaymath}
T_{q_1}\oplus T_{q_2}\rightarrow H^1(D,\mu^*T_\h)\rTo^{\vartheta} 
\bigoplus\limits_{i=1,2} H^1(D_i,\mathcal{L}_{i})\rightarrow 0
\end{displaymath}
\noindent
Since $H^1(D_i,\mathcal{L}_i)$ has dimension $2c_i-1$,
$\vartheta$ is a surjective morphism between two vector spaces of the same
dimension, i.e. it is an isomorphism.
\end{proof}
\begin{proposition}\label{prop:e}
Let $\mathcal{L}_i$ be as in \ref{lem:l} and  
$L_i=\mathcal{E}xt^1(\Omega_D,\mathcal{O}_D)_{q_i}$, $i=1,2$, be the line 
bundle corresponding to the deformations resolving the $i$-th node. 
Then the obstruction bundle $\e$ fits in the exact
sequence:
\begin{equation}\label{eq:seq}
0\rightarrow\bigoplus\limits_{i=1,2}L_i\rightarrow
\bigoplus\limits_{i=1,2} H^1(D_i,\mathcal{L}_{i})\rightarrow\e\rightarrow 0
\end{equation}
\noindent
In particular it has rank $\sum_{i=1}^2(2c_i-2)$.  
\end{proposition}
\begin{proof}
In general the absolute obstruction theory
$\ext^{\bullet}(f^*\Omega_X\to\Omega_{\mathcal{C}},\mathcal{O}_{\mathcal{C}})$
 for $\mbar{g,n}(X,\beta)$ is induced by the relative obstruction theory 
$\ext^{\bullet}(f^*\Omega_X,\mathcal{O}_{\mathcal{C}})$ over 
$\mathfrak{M}_{g,n}$ (see \S \ref{sec:c} for notations).\\
The moduli space $\mbar{0,0}(\h,(1,0,c))$ is not smooth over 
$\mathfrak{M}_{0,0}$; however it is smooth over the smooth $2$-codimensional 
locus defined by not smoothening the nodes at $q_1,q_2$ 
(see Rmk.\ref{rem:liscezza m(c)}). The normal space to this locus is 
$\bigoplus\limits_{i=1,2}\ext^1(\Omega_D,\mathcal{O}_D)_{q_i}$.
(For more details see \cite{p}, Propositions 2.2.4-2.6.10).
\end{proof}
\noindent
For each $i$, let $\e_{c_i}$ be the cokernel of the injection 
$L_i\rightarrow H^1(D_i,\mathcal{L}_{i})$. It is a vector
bundle of rank $2c_i-2$ on $M(c_i)$. It is the one we find when we have
only one node on $D$. Since $\e_{c_1}\oplus\e_{c_2}$ and $\e$ fit into the
same exact sequence, it holds
$c_{\textrm{top}}(\e)=c_{\textrm{top}}(\oplus\e_{c_i})$.\\
In \cite{gr}, Graber constructs a variety $X$ by blowing up $\mathbb{P}^2$ in
a point and then blowing up a point on the exceptional divisor. He gets two
exceptional divisors meeting in a node. Let $A$ be the (-1)-curve, $B$ the
(-2)-curve and $\beta_{c}=A+cB$. He shows that the moduli space
$\mbar{0,0}(X,\beta_{c})$ is smooth of expected dimension zero and isomorphic
to $M(c)$. Besides its virtual fundamental class can be realized as the top
Chern class of a vector bundle $\tilde{\e}_{c}$ which sits in the same exact 
sequence defining the bundle $\e_{c}$. Then $c_{\textrm{top}}(\e_{c})=
c_{\textrm{top}}(\tilde{\e}_{c})$.
\begin{proposition}\emph{\textbf{(Graber)}}\label{prop:vanishing}   
For all $c\geq 2$, $c_{\textrm{top}}(\tilde{\e}_{c})=0$.     
\end{proposition}
\begin{proof}
This is Proposition 3.5 in \cite{gr}.
\end{proof}
\begin{rem}\label{rem:claim}
Let $M^*(c)$ be the fiber over $(0,\infty)$ of the evaluation map
$ev=(ev_1,ev_2):\mbar{0,2}(\mathbb{P}^1,c)\rightarrow\mathbb{P}^1\times
\mathbb{P}^1$. Denote by $\e_c^*$ the obstruction bundle of $M^*(c)$. The
following diagram is commutative:
\begin{displaymath}
\begin{diagram}[height=0.8cm,width=1.0cm]
M^*(c) & \rTo^f & M(c) \\
\dTo & & \dTo \\
\mbar{0,2}(\mathbb{P}^1,c) & \rTo^g & \mbar{0,1}(\mathbb{P}^1,c) 
\end{diagram}
\end{displaymath} 
\noindent
where $g$ and $f$ forget the point mapping to $\infty$. In particular it can
be proved that $\e^*_c$ is the pullback bundle 
$f^*\e_c$ of the obstruction bundle of $M(c)$, so that its top Chern class 
vanishes for $c\geq 2$, (see \cite{p} Rmk 2.6.12, Lem.2.6.13).
\end{rem}
\begin{theorem}\label{thm:classe virtuale}
The virtual fundamental class of a component of a general fiber of $\tau$ 
is given by:
\begin{displaymath}
\begin{array}{ll}
\lbrack M(c_1)\times M(c_2)\rbrack^{vir}=
[M(c_1)\times M(c_2)] & \textrm{if}~~ 0\leq
c_1,c_2\leq 1\\
 & \\
\lbrack M(c_1)\times M(c_2)\rbrack^{vir}=0 & \textrm{otherwise}
\end{array}
\end{displaymath} 
\end{theorem}
\begin{proof}
If $c_1, c_2$ are $0$ or $1$ then $M(c_1)\times M(c_2)$ is smooth of the
expected dimension equal to zero and the virtual fundamental class coincide
with the usual fundamental class. If $c_1$ or $c_2$ is bigger than or equal to
$2$ then by \ref{prop:vanishing} the top Chern class of the obstruction bundle
vanishes. 
\end{proof}
\subsection{Some vanishing results}\label{sec:gwi-1}

We prove some vanishing results for the GW invariants which are related to the
particular geometry of the effective curves involved.\\
For an exhaustive treatment of the invariants and their properties see for 
instance \cite{km}.\\
Note that propositions \ref{prop:obs-bdle} and \ref{rel-obs} imply:   
\begin{theorem}\label{thm:inv}
Let $\pi:\mbar{g,n}(X,\beta)\rightarrow\mbar{g,0}(X,\beta)$ be the usual map
forgetting the markings and $ev=(ev_1,\ldots,ev_n)$ be the evaluation map.
Let $\e$ be the obstruction sheaf on $\mbar{g,0}(X,\beta)$.
Choose cycles $\Gamma_1,\ldots,\Gamma_n$ in $X$ representing the cohomology
classes $\gamma_1,\ldots,\gamma_n$ such that $ev^{-1}_i(\Gamma_i)$ intersect
generically transversally. Then if $A=\pi_*(\cap_i ev^{-1}_i(\Gamma_i))$ is a
cycle in the smooth locus of $\mbar{g,0}(X,\beta)$:
\begin{equation}
\langle \gamma_1\cdot\ldots\cdot\gamma_n\rangle_{\beta}=\int_{A}c_{top}(\e)
\end{equation}
\end{theorem}
\noindent
Consider the classes $T_4=[\Gamma(p)]$, $T_5-T_4$, $\frac{1}{2}T_6$, 
$\frac{1}{2}T_7$ in $\ah{2}$. 
\begin{proposition}\label{prop:van1}
For $\gamma$ equal to one of the above classes, 
$\langle \gamma \rangle_{(0,0,c)}=0$. 
\end{proposition}
\begin{proof}
Suppose $c=1$. A curve $(0,0,1)$ is incident to the cycle $\Gamma(p)$ if it is
the curve of non-reduced subschemes supported on $p$, i.e. if it is the fiber
over $p$ of the support map $s$:
\begin{displaymath}
\D\cong\mbar{0,1}(\h,(0,0,1))\rTo^{s} \mbar{0,0}(\h,(0,0,1))\cong Q
\end{displaymath} 
Let $ev$ be the evaluation map $\mbar{0,1}(\h,(0,0,1))\rightarrow \h$.\\
Since $s$ is flat, $s^*(p)=s^{-1}(p)$ and it is of codimension 2 in
$\mbar{0,1}(\h,(0,0,1))$.\\
As a set $ev^{-1}(\Gamma(p))=s^{-1}(p)$, so $ev^*(T_4)=\lambda s^*(p)$ has
 codimension 2.
\begin{displaymath}
\langle T_4 \rangle_{(0,0,1)}=\int_{[\mbar{0,1}(\h,(0,0,1))]^{vir}}
ev^*(T_4)=\lambda\int_{[\mbar{0,1}(\h,(0,0,1))]}
s^*[p\cdot c_{top}(\e)]=0
\end{displaymath}
where $\e$ is the obstruction bundle on 
$\mbar{0,0}(\h,(0,0,1))$ and  it has rank $1$, 
so that $p\cdot c_{top}(\e)=0$ on $\mbar{0,0}(\h,(0,0,1))$. 
Curves of type $(0,0,c)$ intersecting $\Gamma(p)$ are multiple covers of
$(0,0,1)$, so $\langle T_4\rangle_{(0,0,c)}=0$.\\
The cycle class $T_5-T_4$ can be represented by the set of subschemes whose
support is incident to two lines $l_1, l_2$ with $l_k\in \w{k}$, $k=1,2$.
A curve $(0,0,1)$ can meet such a cycle only if it is the curve supported on
the incident point $l_1\cap l_2$. The previous argument works and
$\langle T_5-T_4\rangle_{(0,0,c)}=0$.\\
The cycle classes $\frac{1}{2}T_6, \frac{1}{2}T_7$ are represented by 
the sets of subschemes 
with support incident to two lines in the same ruling, so a curve $(0,0,1)$ 
can never meet these cycles. This concludes the proof. 
\end{proof}
\begin{lemma}\label{lem:t8}
For each $c\geq 1$, $\langle T_{8}\rangle_{(0,0,c)}=4/c^2$.
Symmetrically, the formula holds also for $T_9$. 
\end{lemma}
\begin{proof}
Consider the diagram:
\begin{displaymath}
\xymatrix{
\mbar{0,1}(\h,(0,0,c))\ar[r]^{\quad \quad \tilde{ev}}\ar[d]_\pi 
& \D\ar[r]^i\ar[d]^s & \h \\
\mbar{0,0}(\h,(0,0,c))\ar[r]_{\quad \quad g} & Q &  }
\end{displaymath}
\noindent
where $g([C,\mu])=\supp{\mu(C)}$. 
Set $ev$ to be the composition map $i\circ\tilde{ev}$. \\
We know that $c_{2c-1}(\e)=-g^*K_Q\cdot c_{2c-2}(\tilde{\e})$, where $\e$ is
the obstruction sheaf on $\mbar{0,0}(\h,(0,0,c))$ and
$\tilde{\e}$ is the sheaf defined in \ref{prop:ec}. So we have to calculate:
\begin{displaymath}
\langle T_{8}\rangle_{(0,0,c)}=\int\limits_{[\mbar{0,1}(\h,(0,0,c))]}
ev^*T_8\cdot\pi^* (g^*(-K_Q)\cdot c_{2c-2}(\tilde{\e}))
\end{displaymath}
\noindent
Note that a point $[C, x, \mu]\in\mbar{0,1}(\h,(0,0,c))$ 
is such that the support of $\mu(C)=\mu(x)=Z$ is a point $p\in Q$, because 
a curve of class $(0,0,c)$ is a multiple cover of a fiber of $s$.\\ 
The above diagram is commutative, let $f$ be the composition
 $g\circ\pi=s\circ\tilde{ev}$. Let $h_1$ be the cycle class of the first
ruling on $Q$ and $\zeta=c_1(\mathcal{N}_{\tilde{\delta}|\tu})$ of degree
$-1$ on a fiber of the blowup map $\tilde{\delta}\rightarrow\delta$.
Then it is easy to verify that 
$i^*T_8=2\cdot s^*h_1\cdot\zeta$. 
We have to calculate the degree:
\begin{displaymath}
\int\limits_{[\mbar{0,1}(\h,(0,0,c))]}
2f^*(-K_Q\cdot h_1)\cdot\tilde{ev}^*\zeta\cdot\pi^* (c_{2c-2}(\tilde{\e}))
\end{displaymath}
\noindent
Since $-K_Q\cdot h_1=2h_3$ where $h_3$ is the point-class in $A^2(Q)$, 
we get:
\begin{displaymath}
f^*(-K_Q\cdot h_1)\cdot\tilde{ev}^*\zeta=2\tilde{ev}^*(\zeta\cdot s^*h_3)
\end{displaymath}
\noindent
Let $x\in\D$ be a point, we denote by $M_1$ the inverse image 
$\tilde{ev}^{-1}(x)$
and by $M_0$ its image $\pi(M_1)=g^{-1}(s(x))$ in $\mbar{0,0}(\h,(0,0,c))$.
The restricted morphism $\tilde{\pi}:M_1\rightarrow M_0$ has degree 
$c$. In particular:
$$\tilde{\pi}_*[M_1]=c[M_0]=c\cdot g^*[s(x)]$$
\noindent 
Since $\zeta\cdot s^*h_3=[x]$ is the point-class in $\D$, by the projection
formula and what we said in \S  \ref{sec:c}, our invariant is:
{\setlength\arraycolsep{2pt}
\begin{eqnarray}
\langle T_{8}\rangle_{(0,0,c)}& = & \int\limits_{[\mbar{0,1}(\h,(0,0,c))]}
4\tilde{ev}^*(\zeta\cdot s^*h_3)\cdot\pi^*c_{2c-2}(\tilde{\e}){}\nonumber \\
   & = & {}\int\limits_{[M_1]}4\tilde{\pi}^*c_{2c-2}(\tilde{\e})
          =\int\limits_{[M_0]} 4c\cdot c_{2c-2}(\tilde{\e}){}\nonumber \\
   & = & {}4c\cdot c_{2c-2}(\tilde{\e}|_{g^{-1}(s(x))})=\frac{4}{c^2}\nonumber
\end{eqnarray}}   
\end{proof}
\noindent
With notations as in section \ref{sec:m(c)}, let $M^*\subseteq
\mbar{0,0}(\h,(1,0,c))$, $c\geq 2$, be the closed subset of stable maps 
$\mu:D\rightarrow\h$ 
where the domain curve is reducible and  
$\mu(D_0)=C(l_1)$ is tangent to the conic defined 
by $\D$ in Sym$^2(l_1)$. Consider the following maps:
\begin{displaymath}
\mbar{0,3}(\h,(1,0,c))\rTo^{\pi}\mbar{0,0}(\h,(1,0,c))\rTo^{\tau}
\mbar{0,0}(\h,(1,0,0))
\end{displaymath}
\noindent
The map $\pi$ forgets the marked points and (eventually) stabilizes the curve. 
The map $\tau$ is the surjective map defined in \S \ref{sec:m(c)}.
It is easy to see that fibers of
$\tau$ are $3$-codimensional (see also Remark 2.8.3 in \cite{p} for a detailed
proof). 
\begin{proposition}\label{prop:van2}
If $c>2$ then all GW invariants 
$\langle\gamma_1\gamma_2\gamma_3\rangle_{\beta}$ 
for curves of type $(1,0,c), (0,1,c)$ vanish.
\end{proposition}
\begin{proof}
The two cases are symmetric. We consider only $(1,0,c)$.\\
We have seen that such a curve is reducible. 
It has a component of class $(1,0,0)$ not contained into 
$\D$ and it decomposes as: 
$$(1,0,c)=(0,0,c_1)+(0,0,c_2)+(1,0,0)$$
\noindent 
with $c_1,c_2\geq 0,~c_1+c_2=c$.\\
We are free to choose a basis of $\ah{*}$ such that every cycle class can be
represented by cycles intersecting the stratification properly. It is enough 
to prove that GW invariants involving such classes vanish. Choose three of
them $\gamma_1, \gamma_2, \gamma_3$ satisfying  $\sum\textrm{cod}~\gamma_i=6$. 
This condition means we are looking at three possible 3-uples 
of elements whose codimensions, up to a permutation of indexes, are 
$(1,1,4),~(1,2,3),~(2,2,2)$.\\
Consider the diagram:
\begin{displaymath}
\begin{diagram}[height=0.8cm,width=1.0cm]
A      & \rTo & \gamma_1\times\gamma_2\times\gamma_3 \\
\dTo   &      & \dTo                    \\
\overline{M}_{0,3}(\h,(1,0,c)) & \rTo^{ev} & \h^3 
\end{diagram}
\end{displaymath}
\noindent
where $A=ev^*(\gamma_1\times \gamma_2\times \gamma_3)$. 
By the Position Lemma cod$~A\geq 6$ (see \ref{rem:weak-pos}).
Let $\pi$ be the flat map defined before and $B=\pi(A)$, then cod$~B
\geq 3.$ If cod$~B>3$, the GW invariants vanishes for dimensional reasons, 
so we can assume cod$~B=3$.\\
If the class of a map $[f]$ is in $B$, then all the points in $\tau^{-1}
(\tau([f]))$ are in $B$, because they differ only by the choice of a multiple
 cover of $(0,0,1)$ and this does not affect incidence conditions. The 
codimension of a fiber of $\tau$ is already equal to 3, so $B$ is a union of 
finitely many components of fibers of $\tau$. 
With notations as in section \ref{sec:m(c)} the set $B$ is:
\begin{displaymath}
B=\coprod\limits_{\stackrel{c_1+c_2=c}{c_i\geq 0}} 
M(c_1)\times M(c_2)
\end{displaymath}
\noindent
where $M(0)$ is a point.
If $c>2$ then there exists $i$ such that $c_i>1$. By \ref{thm:classe virtuale}:
\begin{displaymath}
\langle\gamma_1,\gamma_2,\gamma_3\rangle_{(1,0,c)}=0
\end{displaymath}
\end{proof}
\section{Quantum Cohomology}

Quantum Cohomology is a deformation of the cup product of $\ah{*}$ involving
the genus zero Gromov-Witten Invariants. In particular the Small 
Quantum Cohomology ring $QH_s^*(\h)$ of $\h$ incorporates only the 
genus zero $3$-point GW invariants in its product. \\
In this section we give a presentation of $QH_s^*(\h)$ and describe a
(partial) algorithm computing all the GW invariants on $\h$.\\
We do not give explicitly the complete computations of all the invariants we 
need, for a more detailed treatment we refer to \cite{p}, Chap. 3.\\

\noindent
\textbf{Notations:} the cup product in $\ah{*}$ will be denoted by 
$\alpha\cup\beta$. 
We will use the symbol $\langle T^{n}\rangle_{\mathrm{\beta}}$ to denote the GW
invariant 
$\langle\underbrace{T\cdot\ldots\cdot T}_{n}\rangle_{\mathrm{\beta}}$.
\subsection{The Small Quantum Cohomology Ring}\label{sec:small}

Let $T_0=1,T_1,\ldots,T_{13}$ be a homogeneous $\mathbb{Q}$-basis for $\ah{*}$
such that $T_1,T_2,T_3$ generate $\ah{1}$. We denote by $(g_{ij})$ the 
matrix $(\int_{\h}T_i\cup T_j)$ and by $(g^{ij})$ its inverse. We introduce 
formal variables
$\{y_0,q_1,q_2,q_3,y_4,\ldots,y_{13}\}$ which we will abbreviate as $q,y$. 
For $\beta$ an effective class in $A_1(\h)$, the following expression 
defines a power series in the ring $\mathbb{Q}[[q,y]]$:
\begin{equation}\label{eq:gm}
\gm{ijk}=\sum_{n\geq 0}\sum_{\beta\neq 0}\frac{1}{n!}\langle\gamma^n T_iT_jT_k
\rangle_{\beta}\cdot q^{\beta} 
\end{equation}
\noindent
where $\gamma=y_4T_4+\cdots+y_{13}T_{13}$ and $q^{\beta}=q_1^b\cdot q_2^a\cdot
q_3^c$ for $\beta=(a,b,c)$, (see \cite{gp}). 
Note that if one of the indexes
 $i,j,k$ is zero, then the expression vanishes, because of the condition
 $\beta\neq 0$.\\
Consider the free $\mathbb{Q}[[q,y]]$-module $\ah{*}\otimes_{\mathbb{Q}}
\mathbb{Q}[[q,y]]$ generated by $T_0,\ldots$, $T_{13}$. The so called 
\emph{$*$-product} yields a $\mathbb{Q}[[q,y]]$-algebra structure on it and it
is defined by:
\begin{displaymath}
T_i*T_j=T_i\cup T_j+\sum_{e,f=0}^{13}\gm{ije}g^{ef}T_f
\end{displaymath}
\noindent
It is well known that the $*$-product is commutative,
associative, with unit $T_0$. In particular let
 $\gamma_1,\ldots,\gamma_n$ be cohomology classes on $\h$, $\beta$ in $
A_1(\h)$ the class of an effective curve and $A,B$ sets of indexes. Then the
associativity law reads:
\begin{gather}\label{eq:ass}
\sum   \langle T_i\cdot T_j\cdot T_e\cdot\prod_{a\in A}
       \gamma_a\rangle_{\beta_1}~g^{ef}
       \langle T_k\cdot T_l\cdot T_f\cdot\prod_{b\in B}
       \gamma_b\rangle_{\beta_2}=\nonumber \\
=\sum  \langle T_i\cdot T_k\cdot T_e\cdot\prod_{a\in A}
       \gamma_a\rangle_{\beta_1}~g^{ef}
       \langle T_j\cdot T_l\cdot T_f\cdot\prod_{b\in B}
       \gamma_b\rangle_{\beta_2}    
\end{gather} 
where the sum is over all the possible partitions $A\cup B=[n]$ of $n$
indexes, all possible sums $\beta_1+\beta_2=\beta$ with $\beta_i$ effective
and over $e,f=0,\ldots,13$. \\
\noindent
By definition the Big
Quantum Cohomology ring of $\h$ is the $\mathbb{Q}[[q,y]]$-algebra 
$(\ah{*}\otimes_{\mathbb{Q}}\mathbb{Q}[[q,y]], *)$. \\
\\
\noindent
The Small Quantum Cohomology ring $QH_s^*(\h)$ of $\h$ is defined by setting 
to zero all the formal variables $y_i$. Moreover $QH_s^*(\h)$ is graded by  
deg~$q_1=$deg~$q_2=2$, deg~$q_3=0$ and deg~$T_j=$cod~$T_j$.
Since $q_1,q_2$ have positive degree, we have 
$T_i* T_j\in \ah{*}\otimes_{\mathbb{Q}}\mathbb{Q}[q_1,q_2][[q_3]]$. Hence 
the Small Quantum Cohomology ring of $\h$ is:
\begin{displaymath}
QH_s^*(\h)=(\ah{*}\otimes_{\mathbb{Q}}\mathbb{Q}[q_1,q_2][[q_3]],*)
\end{displaymath}
\noindent
It is a deformation of $\ah{*}$ in the usual sense, in
fact we can recover the Chow ring of $\h$ by setting all the $q_i$ variables 
equal to zero.\\
Let $\mathbb{Q}[Z]=\mathbb{Q}[Z_1,\ldots,Z_4]$ and let
\begin{displaymath}
\ah{*}=\frac{\mathbb{Q}[Z]}{(f_1,\ldots,f_s)}
\end{displaymath}
be a presentation with arbitrary homogeneous generators $f_1,\ldots,f_s$
for the ideal of relations. Finally let $\mathbb{Q}(q,Z)=
\mathbb{Q}[q_1,q_2,Z_1,\ldots,Z_4][[q_3]]$. 
The following proposition is a slightly modified version of \cite{fp} \S10 
Prop.11. 
\begin{proposition}\label{prop:presentation}   
Let $f^{'}_1,\ldots,f^{'}_s$ be any homogeneous elements in 
$\mathbb{Q}(q,Z)$ such that:
\begin{itemize}
\item[(i)] $f^{'}_i(0,0,0,Z_1,\ldots,Z_4)=f_i(Z_{1},\ldots,Z_4)$ in ~ 
$\mathbb{Q}(q,Z)$,
\item[(ii)] $f^{'}_i(q_1,q_2,q_3,Z_1,\ldots,Z_4)=0$ in ~$QH_s^*(\h)$.
\end{itemize} 
\noindent
Then the canonical map 
\begin{displaymath}
\frac{\mathbb{Q}(q,Z)}{(f^{'}_1,\ldots,f^{'}_s)}\rTo^{\varphi} 
QH_s^*(\h)
\end{displaymath}
\noindent
is an isomorphism.
\end{proposition}
\begin{proof}
As in \cite{fp} we can use a Nakayama-type induction. First we observe that
given a homogeneous map $\psi:M\rightarrow N$ between two finitely generated 
$\mathbb{Q}(q,Z)$-modules such that the induced map:
\begin{displaymath}
\frac{M/(q_3)}{(q_1,q_2)}\rTo^{\psi_{1,2}} \frac{N/(q_3)}{(q_1,q_2)}
\end{displaymath}
\noindent is surjective, then $\psi_3:M/(q_3)\rightarrow N/(q_3)$ is 
surjective, because $q_1,q_2$ have positive degree. Since the ideal $(q_3)$ 
is contained into the radical of Jacobson of
$\mathbb{Q}(q,Z)$ and $N=\psi(M)+(q_3)N$, by surjectivity of $\psi_3$, it
follows that $\psi$ is surjective (\cite{am} Cor. 2.7). Hence by 
hypothesis $(i)$ our map $\varphi$ is surjective. If $\tilde{T}_i$,
$i=0,\ldots,13$ are homogeneous lifts to $\mathbb{Q}[q_1,q_2][[q_3]]$ of a 
basis of $\ah{*}$, exactly the same argument of passing to the quotients shows
that their images in $\mathbb{Q}(q,Z)/(f^{'}_1,\ldots,f^{'}_s)$ generates 
this $\mathbb{Q}[q_1,q_2][[q_3]]$-module. But $QH_s^*(\h)$ is free over 
$\mathbb{Q}$ of rank $14$, so $\varphi$ is an isomorphism.
\end{proof}
\subsection{Explicit calculations of some invariants}\label{sec:gwi-2}

We calculate all the monomials arising from the 
$*$-product of two generators of $\ah{*}$, except $T_4*T_4$.\\
We distinguish different cases.
{\setlength\arraycolsep{2pt}
\begin{eqnarray}
T_3*T_3 & = & T_3\cup T_3+\sum\limits_{c\geq 1} 
2\left\{ 2T_5+T_6+T_7-T_8-T_9\right\} q_3^c +{} \nonumber\\
  &  & {}+\sum\limits_{c\geq 0} c^2
\left\{\langle T_{13}\rangle_{(1,0,c)}q_2q_3^c+\langle T_{13}\rangle_{(0,1,c)}
q_1q_3^c\right\} T_0 \nonumber
\end{eqnarray}}
\noindent
where we use \ref{prop:van1},~\ref{lem:t8} and \ref{rem:sym}.\\
If $T_i$ is a divisor class with $i\neq 3$:
\begin{displaymath}
T_i*T_3=T_i\cup T_3+\sum\limits_{c\geq 0}c \langle T_{13}
\rangle_{\beta}\cdot T_0\cdot q^\beta ~~~\textrm{with}~
\beta=(1,0,c)~\textrm{or}~(0,1,c)
\end{displaymath}
\noindent
If $T_i,T_j$ are divisor classes with $i,j\neq 3$:
\begin{displaymath}
T_i*T_j=T_i\cup T_j+\sum\limits_{c\geq 0}\langle T_iT_jT_{13}
\rangle_{\beta}\cdot T_0\cdot q^\beta ~~~\textrm{with}~
\beta=(1,0,c)~\textrm{or}~(0,1,c)
\end{displaymath}
\noindent
If $T_i$ is a divisor class with $i\neq 3$:
\begin{displaymath}
T_i*T_4=T_i\cup T_4+\sum\limits_{\stackrel{\cod T_e=3}{c\geq 0}}
\langle T_iT_4T_e
\rangle_{\beta} g^{ef} T_f\cdot q^\beta ~~~\textrm{with}~
\beta=(1,0,c)~\textrm{or}~(0,1,c)
\end{displaymath}
\noindent
Finally:
\begin{displaymath}
T_3*T_4=T_3\cup T_4+\sum\limits_{\stackrel{\cod T_e=3}{c\geq 0}}
c \langle T_4T_e
\rangle_{\beta}g^{ef} T_f\cdot q^\beta ~~~\textrm{with}~
\beta=(1,0,c)~\textrm{and}~(0,1,c)
\end{displaymath}
\noindent
where we use \ref{prop:van1} again.\\ 
By the vanishing result \ref{prop:van2}, it is enough to calculate 
$\langle T_{13} \rangle_{\beta}$ and $\langle T_{4},\cod 3 \rangle_{\beta}$
with $\beta=(1,0,c),(0,1,c),~0\leq c\leq 2$. 
\begin{theorem}\label{thm:some-inv}
It holds 
$\langle T_{13}\rangle_{(1,0,1)}=2$, $\langle T_{4}T_{10}\rangle_{(1,0,1)}=
\langle T_{4}T_{12}\rangle_{(1,0,1)}=1$. \\
All other invariants of the form 
$\langle T_{13} \rangle_{(1,0,c)}$ and $\langle T_{4},\cod 3 \rangle_{(1,0,c)}$
 are zero.
\end{theorem}
\begin{proof}
We give a detailed proof for the case $\langle T_{13} \rangle_{(1,0,c)}$; the
other proofs are similar, for details see \cite{p}, \S 3.4.\\  
If $\e$ is the rank $\ds{\h}-3$ obstruction bundle on
$\mbar{0,0}(\h,(1,0,c))$, by Theorem \ref{thm:inv} we have to compute:
{\setlength\arraycolsep{2pt}
\begin{eqnarray}
\langle T_{13}\rangle_{(1,0,c)} 
 & = & {}\int_{ev^{-1}(Z)}\pi^*c_{\ds{\h}-3}(\e) {}\nonumber
\end{eqnarray}}
\noindent
where $Z$ is a generic point of $\h$ representing the class $T_{13}$ and 
$\pi$ is the map forgetting a point and stabilizing.\\
If $c=0$, we know that $\mbar{0,0}(\h,(1,0,0))$ is smooth of the 
expected dimension $\ds{\h}=3$.
In particular the top Chern class of $\e$ gives $1$. 
We can choose a representative $Z$ of the class
$T_{13}$ such that $l_Z\notin W_1$, 
then the fiber $ev^{-1}(Z)$ is empty and the
GW invariant vanishes.\\
If $c=1$, we have to analyse separately 
what happens on the two components of the moduli space. We can 
choose $Z\notin\D\cup\mo$, with $\supp{Z}=\{p_0,q_0\}$, 
so that reducible curves of type
$(1,0,1)$ give no contribution to the invariant. Let us consider a
stable map with image an irreducible curve. 
It is a smooth
point for the moduli space $\mbar{0,1}(\h,(1,0,1))$ 
which is $4$-dimensional in it. 
Denote by $M^{irr}$ the irreducible component parametrizing such maps, 
then $ev(\overline{M^{irr}})=\h$. 
The restricted map $ev:\overline{M^{irr}}\rightarrow\h$ has degree two, 
because an irreducible curve $C$ of class $(1,0,1)$ is completely determined 
by choosing a line $l_1\in W_1$ and a point $p_1\notin l_1$ and all its points 
are reduced. Hence the fiber over $Z$ contains
two points: the isomorphism classes $[\mathbb{P}^1, x, \mu]$ where 
$\mathbb{P}^1\cong C(p_0,l_1(q_0))$ or $\mathbb{P}^1\cong C(q_0,l_1(p_0))$.
In the first case, $\mu$ is defined by $\mu(t)=(p_0, f(t))$, 
with $f:\mathbb{P}^1\rightarrow Q$
a parametrization of $l_1(q_0)$ such that $f(x)=q_0$. 
Similarly for the other map.
Then we have a contribution equal to $2$ to the GW invariant.\\
If $c=2$, all the curves of class $(1,0,2)$ are reducible contained into
$\D\cup\mo$, choosing $Z\notin\D\cup\mo$ the fiber $ev^{-1}(Z)$ is
empty and the GW invariant vanishes.
\end{proof}
\noindent  
We want to apply proposition \ref{prop:presentation} to get a presentation of
$QH_s^*(\h)$. Hence we need to write down the $17$ relations defining $\ah{*}$
using the $*$-product. We will denote them by $f_i^*$. 
Using associativity, we can calculate almost all the GW 
invariants we need to reach our aim.\\
For example the identity $(T_1*T_1)*T_2=T_1*(T_1*T_2)$ gives:
\begin{displaymath}
2q_1q_3T_2 +\sum\limits_{\stackrel{\cod T_e=3}{c\geq 0}}
\langle T_6T_e\rangle_{(1,0,c)}g^{ef}T_f\cdot q_2q_3^c=
\sum\limits_{\stackrel{\cod T_e=3}{c\geq 0}}\langle
T_5T_e\rangle_{(0,1,c)}g^{ef}T_f\cdot q_1q_3^c 
\end{displaymath}
By comparing the coefficients of the variables and by \ref{prop:van2} we find:
\begin{displaymath}
\begin{array}{|l|lll|}
\hline
    & & & \\
c=0 & \langle T_5T_e\rangle_{(0,1,0)}=0 & \langle T_6T_e\rangle_{(1,0,0)}=0
& \textrm{for all}~ T_e\in\ah{3}\\
    & & & \\
\hline
    & & & \\
c=1 & \begin{array}{l}
      \langle T_5T_{10}\rangle_{(0,1,1)}=0\\
      \langle T_5T_{11}\rangle_{(0,1,1)}=2\\
      \langle T_5T_{12}\rangle_{(0,1,1)}=2\\
      \end{array}                            &  
      \langle T_6T_e\rangle_{(1,0,1)}=0 & \textrm{for all}~T_e\in\ah{3}\\
    & & & \\
\hline
    & & & \\
c\geq 2 & \langle T_5T_e\rangle_{(0,1,c)}=0 & \langle T_6T_e\rangle_{(1,0,c)}=0
& \textrm{for all}~ T_e\in\ah{3}\\
    & & & \\
\hline                  
\end{array}
\end{displaymath}
\noindent
Among the necessary invariants which can not be computed with this
technique there are those already  calculated in Lemma \ref{lem:t8} and Theorem
\ref{thm:some-inv}. The rest of them can be worked out by hand like in the
following examples. For more details see \cite{p} \S 3.4.
\subsection*{The invariant \emph{$\langle T_{11}T_6\rangle_{(0,1,c)}$}}

We want to calculate:
\begin{displaymath}
\langle T_{11}T_6\rangle_{(0,1,c)}=\int\limits_{[\mbar{0,2}(\h,(0,1,c))]}
ev_1^*(C(p_1,l_1))\cdot ev_2^*(\gamma)\cdot\pi^*(c_{\ds{\h}-3}(\e))   
\end{displaymath}
\noindent
where $\gamma=\{Z\in\h: \supp{Z}\cap l_1^{'}\neq\emptyset, 
\supp{Z}\cap l_1^{''}\neq\emptyset\}$ is a cycle representing $T_6$, for fixed
lines $l_1^{'}, l_1^{''}\in W_1$, and $\e$ is the obstruction bundle on 
$\mbar{0,0}(\h,(0,1,c))$. Both representatives of $T_6$ and $T_{11}$ can be
choosen generic.\\
If $c=0$ the invariant gives $1$, because of the geometry of a curve of 
class $(0,1,0)$.\\ 
If $c=1$ we do not have any contribution from the irreducible curves 
by the genericity
assumptions. Let $C$ be a reducible curve of class $(0,1,1)$. 
It has to be a union $C(l_2)\cup C(p)$ for some $p\in Q$ and $l_2\in
W_{2}$. Since all the points on $C(p_1,l_1)$ and $\gamma$ are reduced, $C$ can
intersects them only along $C(l_2)$. The line $l_2=l_2(p_1)$ 
is then determined. The curve $C(l_2)$ is the line in Sym$^2(l_2(p_1))$
through $(p_1, l_2\cap l_1)$ and $(l_2\cap l_1^{'}, l_2\cap
l_1^{''})$. Moreover there are two possible points for attaching $C(p)$. This
gives a contribution $2$ to the invariant. \\
If $c=2$ we know that each stable map $\mu$ has a reducible domain curve $D$,
in particular $\mu(D)$ is a curve of class $\cdue+ c_1\ctre+
c_2\ctre$ with $c_1+c_2=2$. The $\ctre$-components are points of $M(c_1)$ and
$M(c_2)$ respectively. As before the intersection points with
$C(p_1,l_1)$ and $\gamma$ lie on the $\cdue$-component which is completely
determined. It intersects $\D$ in at most two points $Z_i$, with
$\supp{Z_i}=q_i$. Then by proposition \ref{prop:van2} there is only a point
satisfying all the incident conditions 
$[D, x_1, x_2, \mu]\in\mbar{0,2}(\h,(0,1,2))$:
\begin{displaymath}
\begin{array}{lll}
D=D_0\cup D_1\cup D_2 &  & D_0=C(l_2(p_1))\\ 
\mu_*[D_0]=\cdue &  &\mu_*[D_i]=[C(q_i)],~i=1,2\\
\mu(x_1)=(p_1, l_2(p_1)\cap l_1) & & 
\mu(x_2)=(l_2(p_1)\cap l_1^{'}, l_2(p_1)\cap l_1^{''})\\
y_i=D_0\cap D_i,~i=1,2 & & \mu(y_i)=q_i,~i=1,2
\end{array}
\end{displaymath}
It is a reduced point so it counts with multiplicity one.
\subsection*{The invariant \emph{$\langle T_{13}, \cod 3\rangle_{(1,1,1)}$}}

Choosing generic representatives for the
classes $T_{10}$, $T_{11}$, $T_{12}$, $T_{13}$, stable maps from 
reducible curves of 
class $(1,1,1)$ give no contribution because the expected dimension
of $\mbar{0,2}(\h,(1,1,1))$ is $7$ while reducible curves have less moduli.
Then we restrict to study what happens on the component 
$\mbar{0,2}(\h,(1,1,1))^{irr}$ parametrizing maps from irreducible curves of 
class $(1,1,1)$, which is smooth of the expected dimension. 
Fix a generic point $Z_0$ of $\h$ representing $T_{13}$ with 
Supp$~Z_0=\{p_0,q_0\}$. 
\begin{lemma}\label{lem:claim}
If $(ev_1,ev_2):\mbar{0,2}(\h,(1,1,1))^{irr}\rightarrow\h\times\h$ is 
the evaluation map and $A=\{Z\in\h: l_Z\cap l_{Z_0}\neq\emptyset\}$. Then
$[A]=(ev_2)_*ev_1^*[Z_0]$. 
\end{lemma}
\begin{proof}
Let $(C, x_1, x_2, \mu)\in ev_1^{-1}(Z_0)$ with 
$Z_0=\mu(x_1)$ and $Z_1=\mu(x_2)$. The map $\mu$ is an isomorphism with the 
image curve $\Ll$, which is a line $l$ in
Hilb$^2(\Lambda\cap Q)$ for $\Lambda$ generic plane in $\mathbb{P}^3$. 
Since both 
$Z_0$ and $Z_1$ are in $\Lambda\cap Q$, $l_{Z_1}\cap l_{Z_0}\neq\emptyset$, 
because they lie on the same plane, then $ev_2(ev_1^{-1}Z_0)\subseteq A$. 
The set $ev_1^{-1}(Z_0)$ is
$3$-dimensional as well as $A$, in particular $[A]=T_3$. 
The map $ev_2$ has degree 1 over $A$, in fact given a generic point $Z\in A$,
the lines $l_{Z_0}, l_{Z}$  generate a unique plane $\Lambda$. It
cuts a section $\Lambda\cap Q$ on $Q$ and there is a unique line $l
\subseteq$Hilb$^2(\Lambda\cap Q)$ through $Z_0, Z$. A curve $\Ll$ with two
markings is uniquely determined. Hence the fiber over $Z$
consists of a unique point $[\mathbb{P}^1, x_1, x_2, \mu]$ where
$\mu:\mathbb{P}^1\rightarrow\h$  is an isomorphism with $\Ll$ such that
 $\mu^{-1}(Z_0)=x_1$, $\mu^{-1}(Z)=x_2$. Then $ev_2(ev_1^{-1}(Z_0))$ is
 $3$-dimensional. This proves the lemma.
\end{proof} 
\begin{corollary}\label{lem:7}
For all $T_e\in\ah{3}$:
\begin{displaymath}
\langle T_{13}T_e\rangle_{(1,1,1)}=\int_{\h}T_3\cdot T_e
\end{displaymath}
\end{corollary}
\begin{proof}
It follows from \ref{lem:claim} and \ref{thm:inv}.
\end{proof}
\subsection{A presentation of \emph{$QH_s^*(\h)$}}\label{pres-small}

By the results of the previous section, we can now write down the equations 
$f^*_i$, $i=1,\ldots,17$. Note that we do not write the symmetric equations
obtained by simply interchanging $T_1$ and $T_2$.\\
\\
$T_3*T_3-(T_1+T_2)*T_3+T_1*T_2-\sum\limits_{c\geq 1}
 2q_3^c(2T_1T_2+T_1^2+T_2^2-T_1T_3-T_2T_3)=0$  \\
\\
$T_1*T_1*T_1-q_1(T_3-T_1)+2q_1q_3(2T_1+T_2)+q_1q_3^2(T_1+2T_2-T_3)=0$\\
\\
$T_1*T_1*T_2-2T_1*T_4=0$\\
\\
$T_1*T_1*T_3-2T_1*T_4-2q_1q_3(T_1+T_3)-2q_1q_3^2(T_1+2T_2-T_3)=0$ \\
\\
$T_1*T_2*T_3-2T_3*T_4=0$\\
\\
$T_4*T_4*T_4-2q_1q_2q_3^2T_4=0$\\
\\
$T_4*T_4*T_1-2q_1q_3T_2T_4-q_1q_2q_3T_3+2q_1q_2q_3^2(2T_1+T_2)+$\\
\hspace*{1.9cm} $-q_1q_2q_3^3(2T_1+2T_2-T_3)=0$\\
\\
$T_4*T_4*T_3-2(q_1q_3T_2T_4+q_2q_3T_1T_4)-q_1q_2q_3T_3
            -2q_1q_2q_3^2(2T_1+2T_2-T_3)+$
\hspace*{1.9cm} $-3q_1q_2q_3^3(2T_1+2T_2-T_3)=0$\\
\\
$T_1*T_1*T_4-\frac{1}{2}q_1(T_2T_3-T_1T_2)-q_1q_3(2T_1T_2+T_2^2)
            -\frac{1}{2}q_1q_3^2(T_1T_2+2T_2^2-T_2T_3)+$\\
\hspace*{1.85cm} $-q_1q_2q_3(1+2q_3)T_0=0$\\
\\
$T_1*T_2*T_4-T_4*T_4-q_1q_3T_2^2-q_2q_3T_1^2-q_1q_2q_3(1+2q_3)T_0=0$\\
\\
$T_1*T_3*T_4-T_4*T_4-q_1q_3(T_1T_2+T_2^2+T_2T_3)-q_2q_3T_1^2-
q_1q_3^2(T_1T_2+2T_2^2-T_2T_3)+$\\ 
\hspace*{1.9cm} $-q_1q_2q_3(1+4q_3+3q_3^2)T_0=0$\\
\\
\begin{theorem}\label{thm:small}
The Small Quantum Cohomology ring of $\h$ is: 
\begin{displaymath}
QH^*_s(\h)=\frac{\mathbb{Q}[q_1,q_2,T_1,T_2,T_3,T_4][[q_3]]}
{(f_i^*)_{i=1,\cdots,17}}
\end{displaymath}
\end{theorem}
\begin{proof}
The equations $f^*_i$ satisfy the hypotheses of \ref{prop:presentation}.
\end{proof}
\begin{rem}\label{rem:t4t4}
In the ring $QH^*_s(\h)$ the identity $T_4^2=T_{13}$ 
corresponds to:
$$T_4*T_4=T_{13}+2q_1q_2q_3^2T_0$$
\end{rem}
\subsection{The First Reconstruction Theorem}\label{sec:subring}

All the classes in the fixed basis of $\ah{*}$ can be written as some product 
of the divisor classes except $T_4$. Hence if we restrict to the subalgebra
$\s$  of $\ah{*}$ generated by the divisor classes $T_1,T_2,T_3$, we can apply
the First Reconstruction Theorem (FRT) (\cite{km} Theorem 3.1). 
It says we can compute all the genus zero GW invariants
with arguments in $\s$ by knowing few initial values 
corresponding to the invariants of the form:
\begin{displaymath}
\int_{[\overline{M}_{0,2}(\h,\beta)]^{vir}} 
ev^*(\gamma_1\times\gamma_2)
\end{displaymath}
\noindent
with $ev:\overline{M}_{0,2}(\h,\beta)\rightarrow \h^2$ the usual evaluation
map and $\gamma_1,\gamma_2 \in \s$.
Since $\cod ev^*(\gamma_1\times\gamma_2)$ has to be equal to $2a+2b+3$ and 
$\cod\gamma_i\leq 4$ for $i=1,2$, we find the upper-bound $a+b\leq 2$. 
We have only the following cases:
\begin{displaymath}
\begin{array}{c|c|c|c|c|c|c|}
 \beta   & (0,0,c) & (1,0,c) & (0,1,c) & (1,1,c) & (2,0,c) & (0,2,c)\\
\hline 
         & & & & & &  \\
(\cod\gamma_1, \cod\gamma_2)& (1,2) & (1,4) & (1,4) & (3,4) & 
          (3,4) & (3,4) \\
         & & & & & &  \\
    & & (2,3) & (2,3) & & &  \\
\hline
\end{array}
\end{displaymath}
\noindent
In sections \ref{sec:gwi-1} and \ref{sec:gwi-2} we calculated some of these
invariants. The left ones are obtained by means of 
the associativity. Then we know all of
them. This implies that we can calculate all the GW invariants on $\h$ 
without $T_4$ among the arguments.
\subsection{An algorithm for the tree level GW invariants}\label{sec:alg}

\begin{theorem}
Assume we know all the $2$-point invariants
$\langle\gamma_1\gamma_2\rangle_{\beta}$ with $\gamma_1,\gamma_2\in\s$, and
those of the form $\langle T_4^m\rangle_{\beta}$, $m\geq 1$. 
Then we can compute
recursively all the invariants of type 
$\langle T_4^m\gamma_1\cdot\ldots\cdot\gamma_n\rangle_{\beta}$, with
$\gamma_i\in\s$ such that $4\geq\dg\gamma_1\geq\ldots\geq\dg\gamma_n\geq 2$.
\end{theorem}
\begin{proof} 
We use equation (\ref{eq:ass}) and by induction we suppose to know all
the invariants:
\begin{displaymath}
\begin{array}{ll}
\langle T_4^r\gamma_1\cdot\ldots\cdot\gamma_n\rangle_{\beta} & 
\textrm{with}~r<m\\
 & \\
\langle T_4^r\gamma_1\cdot\ldots\cdot\gamma_s\rangle_{\beta} & 
\textrm{with}~r+s<m+n\\
 & \\
\langle T_4^m\tilde\gamma_1\cdot\ldots\cdot\tilde\gamma_n\rangle_{\beta} &
\textrm{with}~ \dg\tilde\gamma_n<\dg\gamma_n\\
& \\
\langle T_4^m\gamma_1\cdot\ldots\cdot\gamma_n\rangle_{\beta^{'}} &
\textrm{with}~ \beta-\beta^{'}>0~\textrm{effective}\\
\end{array}
\end{displaymath}     
\noindent
If $m=0$, there is no problem because each $\gamma_i$ is in $\s$.\\
If $m\geq 1$ and $n=0$, then we know the values by hypothesis.\\
If $m=1$ and $n=1$, then $\gamma_1$ lives necessarily in codimension $3$ 
and we have
already calculated all the invariants in section \ref{sec:gwi-2}.\\
If $m=1$ and $n\geq 2$, we use (\ref{eq:ass}):
\begin{gather*}
\sum   \langle T_i\cdot T_j\cdot T_e\cdot\prod_{a\in A}
       \gamma_a\rangle_{\beta_1}~g^{ef}
       \langle T_k\cdot T_l\cdot T_f\cdot\prod_{b\in B}
       \gamma_b\rangle_{\beta_2}=\nonumber \\
=\sum  \langle T_i\cdot T_k\cdot T_e\cdot\prod_{a\in A}
       \gamma_a\rangle_{\beta_1}~g^{ef}
       \langle T_j\cdot T_l\cdot T_f\cdot\prod_{b\in B}
       \gamma_b\rangle_{\beta_2}    
\end{gather*} 
\noindent
By induction, we know all the invariants with $\beta_i\neq
0$, $i=1,2$. We look only to the terms with either 
$\beta_1$ or $\beta_2$ equal to zero, i.e. on the left-hand side:
\begin{displaymath}
\underbrace{
\langle T_i\cdot T_j\cdot T_k\cup T_l\cdot
\prod\limits_{1}^n\gamma_s\rangle_{\beta}}_{I_1}+ 
\underbrace{
\langle T_i\cup T_j\cdot T_k\cdot T_l \cdot
\prod\limits_1^n\gamma_s\rangle_{\beta}}_{I_2}
\end{displaymath}
\noindent
on the right-hand side:
\begin{displaymath}
\underbrace{
\langle T_i\cdot T_k\cdot T_j\cup T_l \cdot
\prod\limits_1^n\gamma_s\rangle_{\beta}}_{I_3}+
\underbrace{
\langle T_i\cup T_k\cdot T_j\cdot T_l \cdot
\prod\limits_1^n\gamma_s\rangle_{\beta}}_{I_4}
\end{displaymath}    
\noindent
Since $\gamma_i\in\s$, there exists a decomposition 
$\gamma_n=\alpha\cup\alpha_1$ with
$\alpha_1\in\ah{1}$ and $\dg\alpha=\dg\gamma_n-1$. We choose:
\begin{displaymath} 
\begin{array}{lllll}
T_i=T_4, & T_j=\gamma_1, & T_k=\alpha, & T_l=\alpha_1, &
R=\gamma_2\cdot\ldots\cdot\gamma_{n-1}\\ 
\end{array}
\end{displaymath}
\noindent
Then $I_1$ is the value  
$\langle T_4\gamma_1\cdot\ldots\cdot\gamma_n\rangle_{\beta}$ we want to know 
(this will always be the case). Up to a scalar 
(possibly zero) $I_2$ is 
$\langle T_4\cup\gamma_1\cdot\alpha\cdot R\rangle_{\beta}$, all its arguments
are in $\s$. Analogously $I_4$ is proportional to the known invariant 
$\langle T_4\cup\alpha\cdot\gamma_1\cdot R\rangle_{\beta}$. Finally in
$I_3=\langle T_4\cdot\alpha\cdot\gamma_1\cup\alpha_1\cdot R\rangle_{\beta}$ 
the minimal degree decreased by one. Then we can write $I_1$ as a combination 
of lower degree terms. After a finite number of steps we can reduce our
problem to the previous case with $n=1$.\\ 
If $m\geq 2$ and $n=1$, then we have three possibilities for $\cod\gamma_1$.
If $\cod\gamma_1=4$, we can suppose $\gamma_1=T_{13}$. We choose:
\begin{displaymath} 
\begin{array}{l}
T_k,T_l\in\ah{2}\cap\s~\textrm{with}~T_k\cup T_l=T_{13} \\
T_i=T_j=T_4  \\ 
R=T_4^{m-2} \\
\end{array}
\end{displaymath}
\noindent
We obtain that in $I_2=\langle T_4^{m-2}T_{13}T_kT_l\rangle_{\beta}$ we have a
lower number of $T_4$'s as well as in $I_3$ and $I_4$, since $T_4\cup T_k$,
$T_4\cup T_l$ are in $\s$. We can reduce the problem to find 
$\langle T_4T_{13}\gamma\rangle_{\beta}$, with $\gamma\in\ah{2}\cap\s$, i.e. 
$m=1$.\\
If $\cod\gamma_1=3$, then we can decompose it as
$\gamma_1=\alpha\cup\alpha_1$, with
$\alpha_1\in\ah{1}$ as above. Fixing:
\begin{displaymath} 
\begin{array}{lllll}
T_i=T_4, & T_j=T_4, & T_k=\alpha, & T_l=\alpha_1, &
R=T_4^{m-2}\\ 
\end{array}
\end{displaymath}
\noindent 
we get $I_2$ proportional to $\langle T_4^{m-2}T_{13}\alpha\rangle_{\beta}$,
and we know it by induction. The invariant $I_3=\langle T_4^{m-1}\alpha
\cdot T_4\cup\alpha_1\rangle_{\beta}$ has less $T_4$-classes and the minimal
degree is lower. Finally $I_4$ is proportional to $\langle
T_4^{m-1}T_{13}\rangle_{\beta}$, then it is known.\\
If $\cod\gamma_1=2$, we use the same trick with:
\begin{displaymath} 
\begin{array}{lllll}
T_i=T_4, & T_j=T_4, & T_k=\alpha_1, & T_l=\alpha_2, &
R=T_4^{m-2}\\ 
\end{array}
\end{displaymath}
\noindent   
where $\alpha_1,\alpha_2$ are two divisors such that
$\alpha_1\cup\alpha_2=\gamma_1$. Also in this case we can reduce our problem
to the case $m=1$.\\
If $m\geq 2$ and $n\geq 2$, then we write $\gamma_n=\alpha\cup\alpha_1$,
$\alpha_1\in\ah{1}$ and we choose:
\begin{displaymath} 
\begin{array}{lllll}
T_i=T_4, & T_j=\gamma_1, & T_k=\alpha, & T_l=\alpha_1, &
R=T_4^{m-1}\gamma_2\cdot\ldots\cdot\gamma_{n-1}\\ 
\end{array}
\end{displaymath}
\noindent   
Then $I_2,I_4$ are invariants with less $T_4$'s and in $I_3$ the minimal
degree is $\dg\alpha=\dg\gamma_n-1$. By induction we
reduce to the case $n=1$ or $m=1$.
\end{proof}
\begin{rem}\label{rem:m-odd}
For dimensional reasons, the invariant $\langle T_4^m\rangle_{\beta}$
vanishes unless $m$ is odd. Moreover we know that for $m=1,3$ it is zero.
\end{rem}
\section{Enumerative applications}

We use the results on the Small Quantum Cohomology obtained in the previous
section to count how many hyperelliptic curves on $Q$ of given genus and 
bi-degree  pass through a fixed number of generic points. Basically we reduce a
question in higher genus to a question about rational curves on the Hilbert
scheme $\h$, as in \cite{gr}. To do this we need a relationship between our 
hyperelliptic curves and some rational curves on $\h$.\\  
By hyperelliptic curve we mean a smooth irreducible projective
curve with a choice of hyperelliptic involution, i.e. one with rational
quotient. This involution is unique if the genus is greater than or equal to
$2$.
\subsection{The moduli space of hyperelliptic curves mapping to $Q$}

We recall Lemma 2.1 from \cite{gr}.
\begin{lemma}\label{lem:graber}
If $f:C\rightarrow \mathbb{P}^r$ is a morphism from a hyperelliptic curve such
 that it does not factor through the hyperelliptic map $\pi:C\rightarrow
\mathbb{P}^1$  then $H^i(C,f^*\mathcal{O}(1))$ vanishes for all $i>0$.
\end{lemma}
\noindent
A similar result holds for maps to $Q$.
\begin{lemma}\label{lem:cohom}
Let $p_i:Q\rightarrow\mathbb{P}^1$ be the two projections and 
$\mu:C\rightarrow Q$ be a morphism from a hyperelliptic curve such that 
 $\mu_i:=p_i\circ\mu:C\rightarrow \mathbb{P}^1$, $i=1,2$, 
does not factor through the 
hyperelliptic map. \\
Then $H^i(C,\mu^*T_Q)=0$ for all $i>0$.
\end{lemma}
\begin{proof}
Consider the Euler sequence:
\begin{displaymath}
0\rightarrow \mathcal{O}\rightarrow\mathcal{O}(1)^{\oplus2}\rightarrow 
T_{\mathbb{P}^1}\rightarrow 0
\end{displaymath}
\noindent
Since $T_Q=p_1^*(T_{\mathbb{P}^1})\oplus p_2^*(T_{\mathbb{P}^1})$, 
a surjection is defined:
\begin{displaymath}
H^1(C,\mu^*p_1^*\mathcal{O}^{\oplus2}(1)\oplus \mu^*p_2^*
\mathcal{O}^{\oplus2}(1))\rightarrow H^1(C,\mu^* T_Q)\rightarrow 0
\end{displaymath}
By hypothesis  $H^j(C,\mu^*p_i^*\mathcal{O}^{\oplus2}(1))=0$ for $j>0$,
 so  $H^j(C,\mu^* T_Q)=0$.
\end{proof}
\noindent
Let $M_{g,0}(Q,(d_1,d_2))$ be the moduli space of maps
$\mu:C\rightarrow Q$ from a smooth irreducible projective curve $C$ of 
genus $g$ such that
$\mu_*[C]=(d_1,d_2)$. Let $M_g$ be the moduli space of
semistable projective curves of genus $g$. 
We denote by $H_g$ the sub-locus parametrizing hyperelliptic curves. 
If $C$ is
hyperelliptic then the cyclic group of order $2$ acts on the space of
universal deformations $\mathcal{U}$ of $C$. 
It can be proved that the fixed locus $V\subseteq\mathcal{U}$ is the universal
deformation space of $C$ as a hyperelliptic curve and it is obviously smooth.
It follows that $H_g\subseteq M_g$ is a smooth substack.
The cartesian diagram:
\begin{diagram}[height=0.8cm,width=1.0cm]
\tilde{H}_g(Q,(d_1,d_2))& \rTo & H_g\\
\dTo &                   &\dTo \\
M_{g}(Q,(d_1,d_2)) &\rTo & M_g
\end{diagram}
\noindent
defines the space $\tilde{H}_g(Q,(d_1,d_2))$ parametrizing maps
$\mu:C\rightarrow Q$ from a hyperelliptic curve $C$ of genus $g$ with
$\mu_*[C]=(d_1,d_2)$. We are interested in the open subset $H_g(Q,(d_1,d_2))$ 
of maps $\mu$ such that the composition maps 
$\mu_i=p_i\circ\mu:C\rightarrow \mathbb{P}^1,~i=1,2$ 
do not factor through the hyperelliptic map. 
\begin{theorem}
The natural morphism $\nu:H_g(Q,(d_1,d_2))\rightarrow H_g$ is smooth.
\end{theorem}
\begin{proof}
It follows from the vanishing result \ref{lem:cohom}; for each 
$\mu: C\rightarrow Q$ in $H_g(Q,(d_1,d_2))$, we have $H^1(C,\mu^*T_Q)=0$. 
Then by theorem \ref{thm:smooth}, the forgetful morphism 
$\overline{M}_{g,0}(Q,(d_1,d_2))\rightarrow\mathfrak{M}_{g}$ 
is smooth in
$[\mu]$. Since smoothness is a local property, the
theorem follows. 
\end{proof}
\begin{corollary}
$H_g(Q,(d_1,d_2))$ is smooth and irreducible.
\end{corollary}
\begin{proof}
Smoothness is a direct consequence of the theorem, since both $H_g$
and $\nu$ are smooth.\\
Since $H_g$ is irreducible, it is enough to prove the fibers 
of $\nu$ are irreducible of constant dimension. 
A fiber $\nu^{-1}(C)$ is the set of all $\mu:C\rightarrow Q$ of
bi-degree $(d_1,d_2)$
such that both $\mu_1,\mu_2$ do not factor through the hyperelliptic
map. They are two morphisms to the projective line, so they correspond to two
line bundles on $C$ of degree $d_1,d_2$ respectively. 
We get a morphism
 $f=(f_1,f_2):\nu^{-1}(C)\rightarrow$ Pic$^{d_1}(C)\times$Pic$^{d_2}(C)$.
By Lemma \ref{lem:graber} $Im(f_i)$ is a subset of $\{\mathcal{L}_i:
\mathcal{L}_i~\textrm{is spanned},~h^1(\mathcal{L}_i)=0\}$. Conversely,
for $i=1,2$, let  $W_i$ be the subset of Pic$^{d_i}(C)$ of sheaves 
$\mathcal{L}_i$ such that $\mathcal{L}_i$ is spanned, $h^1(\mathcal{L}_i)=0$
and $\mathcal{L}_i$ is not a multiple of $g^1_2$. Then each $\mathcal{L}_i\in
W_i$ is in the image $Im(f_i)$. $W_i$ is open and dense (if not empty), because
Pic$^{d_i}(C)$ is irreducible. Hence $Im(f_i)$ contains the open subset $W_i$
and therefore it is irreducible (because $W_i$ is). It follows that $Im(f)$ is
irreducible of dimension $2g$. Each fiber 
$f^{-1}(\mathcal{L}_1,\mathcal{L}_2)$,
$\mathcal{L}_i\in W_i$, is a product
$V_1\times V_2$, where $V_i$ is the open set of pairs of global 
sections $(s_i^1,s_i^2)$ of $\mathcal{L}_i$ without common zeros, modulo
scalars. Hence these fibers are irreducible and they have the same
dimension equal to $2(d_1+d_2)-2g$, because the first cohomology of 
$\mathcal{L}_i$ vanishes. Therefore $\nu^{-1}(C)$ is irreducible of 
dimension $2(d_1+d_2)$.
\end{proof}
\subsection{The basic correspondence}\label{sec:corr}

Let $g:\mathbb{P}^1\rightarrow\h$ be a map in $M_{0,0}(\h,(a,b,c))$ 
satisfying the following conditions (\dag):
\begin{flushleft}
- $g(\mathbb{P}^1)$ intersects $\D$ transversally\\
- $g(\mathbb{P}^1)$ is not contained in $\mo$\\ 
- $g(\mathbb{P}^1)$ is disjoint from $\D_2$\\
\end{flushleft}
\noindent
We can associate to $g$ a map $\mu:C\rightarrow Q$ by:
\begin{displaymath}
\xymatrix{
C\ar[r] \ar[d]_{\pi} 
\ar@/^/[rr]^{\mu} & \mathcal{U}\ar[d]^{u}_{2:1}\ar[r]&  Q\\
\mathbb{P}^1\ar[r]_{g} & \h  &  }
\end{displaymath}
\noindent
where $\mathcal{U}$ is the universal family. Then $C$ is a 
smooth hyperelliptic curve of genus $g_{C}=a+b-c-1$ and  
 bi-degree $(b,a)$.
The map $\mu:C\rightarrow Q$ satisfies the following conditions (\ddag):
\begin{flushleft}
- $C$ is a smooth hyperelliptic curve \\
- both $\mu_i$ do not factor through $\pi$\\
- both differentials $d\mu_i$ are injective on ramification points of $\pi$ \\
\end{flushleft}
\noindent
Conversely, let $\mu:C\rightarrow Q$ be an element in $H_g(Q,(d_1,d_2))$, 
with hyperelliptic map $\pi:C\rightarrow\mathbb{P}^1$. If it 
satisfies (\ddag), there exists a
canonical map $g:\mathbb{P}^1\rightarrow \h$ which induces it (see \cite{p},
\S 4.2). More precisely, let us consider the open subset
$M_{0,0}^{\textrm{tr}}(\h,\beta)\subseteq M_{0,0}(\h,\beta) $ of maps from 
irreducible rational curves fulfilling (\dag).
Let $H_g^{\textrm{tr}}(Q,(d_1,d_2))\subseteq H_g(Q,(d_1,d_2))$ be the open 
subset parametrizing maps $\mu$ satisfying (\ddag).
\begin{theorem}\label{theorem:iso}
There is a canonical isomorphism:
\begin{displaymath}
H_g^{\textrm{tr}}(Q,(d_1,d_2))\cong
M_{0,0}^{\textrm{tr}}(\h,(d_2,d_1,d_1+d_2-g-1)) 
\end{displaymath}
\end{theorem}
\begin{proof}
The proof of Theorem 2.4 in \cite{gr} never makes use of the fact that the
curves are in $\mathbb{P}^2$, then it works also for hyperelliptic curves on
$Q$.  
\end{proof}
\subsection{Enumerative results}

By Theorem \ref{theorem:iso} we might expect a relationship between the number
$n$ of hyperelliptic curves on $Q$ of 
bi-degree $(d_1,d_2)$ and genus $g$ 
passing through $r$ general points and some Gromov-Witten 
invariants involving the cycle $\Gamma(p)$, $p\in Q$, and the moduli 
space $\mbar{0,r}(\h,(d_2,d_1,d_1+d_2-g-1))$. 
In particular we need to exclude undesired contributions to the number $n$
 coming from stable maps either living in the wrong dimension or 
with a reducible domain curve.  
\begin{theorem}\label{main}
Fix an effective class $\beta=(a,b,c)\in A_1(\h)$, $a+b\geq 1$, and
 $r$ general points $p_1,\ldots, p_r$ on $Q$ with
$r=2a+2b+1$. Then:
\begin{enumerate}
\item there exists at most a finite number of irreducible rational curves 
of class $\beta$ incident to all the cycles $\Gamma(p_i)$;
\item all such curves intersect $\D\cup\mo$ in points disjoint 
from the $\Gamma(p_i)$;
\item given any arbitrary stable map $\mu:C\rightarrow \h$ of class $\beta$
 incident to all the cycles $\Gamma(p_i)$, then $C$ has a unique irreducible 
component which is not entirely mapped into $\D\cup\mo$, such a component is 
of class $(a,b,c_0)$, where $c_0\leq c$.
\end{enumerate}
\end{theorem}
\begin{corollary}\label{cor:comb}
In the same assumptions, let us consider the usual evaluation map
$ev:\mbar{0,r}(\h,(a,b,c))\rightarrow\h^r$. 
Then $ev^{-1}(\prod\Gamma(p_i))$ is zero dimensional and smooth. 
\end{corollary}
\begin{proof}
Theorem \ref{main} (which is proven in \S \ref{sec:proof}) 
says that given a stable map 
$\mu:C\rightarrow \h$ satisfying all incident conditions, aside from the 
distinguished component of $C$ of
class $(a,b,c_0)$, all other components are of type $(0,0,c')$ and they are 
entirely mapped into $\D$. So they are multiple 
covers of $\mathbb{P}^1$. Moreover, adding a component of type $(0,0,c')$ 
to a stable map can never cause it to be incident to any extra $\G{q}$, 
since it would force another component of the curve to meet the corresponding 
cycle. Finally, different $(0,0,c')$-components are disjoint, since they are
different fibers of the support map $s$, hence they must be incident to the
distinguished component, $C$ been connected. \\
We conclude that the source curve looks like a comb, with the component of
class $(a,b,c_0)$ as the handle and the components of class $(0,0,c')$
as the teeth. We get exactly the same picture obtained in \cite{gr}.\\
There is a finite number of such curves. Infact, if $C$ is irreducible, then 
Theorem \ref{main} confirms our assertion. 
If $C$ is reducible, we have only a finite 
number of possibilities for the multiple covers of a $(0,0,1)$-curve and 
only a finite number of points of intersection of the distinguished component 
with $\D$. So there are only finitely many potential image curves for stable 
maps incident to all of the cycles.\\
To prove smoothness, by Position Lemma \ref{lem:pos}-2 it is enough to 
prove that $ev^{-1}(\Gamma)=ev^{-1}(\Gamma_{reg})$, where 
$\Gamma=\prod\Gamma(p_i)$. \\
By Lemma \ref{lem:trans} $\Gamma_{reg}=
\prod(\Gamma(p_i)-\D)$ and Theorem \ref{main}-2 ensures that 
$ev^{-1}\cap ev_i^{-1}(\Gamma(p_i)\cap\D)$.
\end{proof}
\noindent
With notations as in \S \ref{sec:m(c)}, let $A_1$ be the product
$M(c_1)\times\cdots\times M(c_m)$ with $\sum c_i=c$ and $c_i=1$ for all $i$,
and let $A_2$ be the disjoint union of all the other products with 
$\sum c_i=c$.
\begin{corollary}\label{cor:a}
In the same assumptions, we denote by $A$ the image of 
$ev^{-1}(\prod\G{p_i})$ in $\mbar{0,0}(\h,(a,b,c))$, 
via the usual map forgetting the markings (and stabilizing).  
Then $A=A_1\sqcup A_2$. \\
\noindent
In particular $\dg A_1=\langle T_4^r\rangle_{(a,b,c)}$.
\end{corollary}
\begin{proof}
Theorem \ref{main} says that the only moduli in the choice of a 
stable map
meeting all the $\G{p_i}$ comes from the choice of multiple covers of the 
$(0,0,1)$ curve. Then as a set, each component of $A$ decomposes as a product:
\begin{displaymath} 
M(c_1)\times M(c_2)\times\ldots\times M(c_m)
\end{displaymath}
\noindent
with $c_1+\ldots+c_m=c$. In particular $A$ is contained in the smooth locus of 
$\mbar{0,0}(\h,(a,b,c))$.
Theorems \ref{thm:classe virtuale} and \ref{thm:inv} conclude the proof.  
\end{proof}
\begin{definition}\label{def:comblike}
A curve $C$ in $\h$ of class $(a,b,c)$ is \emph{comblike} if it is 
 a union of an irreducible 
$(a,b,c_0)$-curve, $c_0\leq c$, and 
$c-c_0$ disjoint rational curves, mapping isomorphically onto a 
$(0,0,1)$-curve.  
\end{definition}
\begin{rem}
All the stable maps in $A_1$ have comblike curves as domain.
\end{rem}
\begin{definition}\label{def:enumerative}
Fix $k$ general points $p_i$ on $Q$ 
and $l\geq 0$ general pairs of points $q_j,q'_j$  with $k+3l=r$, 
$r=2d_1+2d_2+1$. 
Let $E^l((d_1,d_2),g)$ be the number of 
hyperelliptic curves on $Q$ of genus $g$ and bi-degree $(d_1,d_2)$ passing 
through all the points and satisfying also the condition that 
$q_i$ is hyperelliptically conjugate to $q'_i$
 for all $i$.
\end{definition}
\begin{theorem}\label{hyp2}
With $\beta=(d_2,d_1,d_1+d_2-g-1)$ and $r,l$ as above:
\begin{equation}\label{eq:formula2}
\langle T_{13}^l\cdot T_4^{r-3l}\rangle_{\beta}=\sum_{h\geq g} {2h+2 
\choose h-g} E^l((d_1,d_2),h)
\end{equation}
\end{theorem}
\begin{proof}
We write $\beta=(a,b,c)$ where $a=d_2, b=d_1, c=d_1+d_2-g-1$.\\
Consider the case $l=0$ and
fix $r$ general points $p_1,\ldots,p_r$. Then the invariant $\langle
T_4^r\rangle_{\beta}$ is given by the degree deg$(ev^{-1}(\prod \Gamma(p_i))$,
 where by Theorem \ref{main} the scheme $ev^{-1}(\prod \Gamma(p_i))$ is
supported in a finite set of points. By Corollary \ref{cor:a}, the only 
contribution to the invariant comes from 
the component of the moduli space
$\mbar{0,r}(\h,\beta)$ corresponding to stable maps from comblike curves such
that the irreducible $(a,b,c_0)$-component is incident to all the cycles
$\Gamma(p_i)$.  
Hence the number of
stable maps is equal to the number of possible irreducible curves of class
$(a,b,c_0)$ times the number of choices for the attachment points of the
$(0,0,1)$-curves. We have to choose $c-c_0$ points among the $2(a+b-c_0)$ ones
in the intersection $(a,b,c_0)\cdot\D$. The formula then follows from the
relationship between $(a,b,c)$ and $(d_1,d_2,g)$.\\
If $l\geq 1$, then the curves have to meet also $l$ general points of $\h$.
Choosing representatives of the point class $T_{13}$ outside $\D\cup\mo$,
curves moving in excess dimension cannot satisfy this condition, so Theorem 
\ref{main} applies also for $l\geq 1$. The same arguments conclude the proof.
\end{proof} 
\begin{rem}\label{rem:finite}
We note that the sum (\ref{eq:formula2}) is finite, in fact the values of $h$ 
are equal to $d_1+d_2-c_0-1$ with $c_0\leq c$. 
By what we showed in section
\ref{sec:alg} if $l\geq 1$ then we can compute all the invariants 
$\langle T_{13}^l\cdot T_4^{r-3l}\rangle_{\beta}$. 
Therefore we can invert the formula
(\ref{eq:formula2}) to get the numbers $E^l((d_1,d_2),h)$.\\
The numbers $E^0((d_1,d_2),g)$ are zero for small values of 
$d_1, d_2, g$. In fact $E^0((d_1,d_2),g)$ is less then or equal to 
$S((d_1,d_1),g)$, the number of
smooth curves of bi-degree $(d_1,d_2)$ of genus $g$ passing through $r$
points. We know that
$S((d_1,d_1),g)$ is zero if $d_1d_2-d_1-d_2-1<0$, hence the first 
possibly nonzero GW invariants with $l=0$ 
are $\langle T_4^{11}\rangle_{(3,2,2)}$ and 
$\langle T_4^{11}\rangle_{(2,3,2)}$.  
\end{rem}
\subsection{Proof of Theorem \ref{main}}\label{sec:proof}

\begin{lemma}\label{lem:disjoint}
With notations as in Theorem \ref{main}, 
let $C$ be an irreducible rational curve
meeting all the cycles $\G{p_i}$ and the orbit $\mo_4$. Then it 
intersects $\D\cup\mo$ in points disjoint from all the $\G{p_i}$.  
\end{lemma}
\begin{proof}
Let $r=2a+2b+1$ and $M\subseteq\overline{M}_{0,r}(\h,(a,b,c))$ be 
the open subset of points $[C,\mu,x_j]_{\scriptstyle{j=1,\ldots,r}}$ such that 
$C\cong\mathbb{P}^1$,
$\mu(C)\cap\mo_4\neq\emptyset$. It is smooth of dimension $2r$. 
The map $M\rightarrow\overline{M}_{0,0}(\h,(a,b,c))$ which
forgets the markings and stabilizes factors through:
\begin{displaymath}
M\rTo^{\pi_i}\overline{M}_{0,1}(\h,(a,b,c))\rTo
\overline{M}_{0,0}(\h,(a,b,c))   
\end{displaymath}
where $\pi_i$ is the map forgetting all the markings but $x_i$ and
stabilizing. It is surjective onto its image Im$(\pi_i)=\mathcal{U}_1$ which is
the universal curve over the smooth locus $\mathcal{U}_0$ of
$\overline{M}_{0,0}(\h,(a,b,c))$. Then $\pi_i:M\rightarrow\mathcal{U}_1$ is 
flat
of relative dimension $r-1$. The set $N=\{[C,\mu,x]: \mu(x)\in\D\cup\mo\}$
is a closed subset of $\mathcal{U}_1$, as it is the inverse image
$ev^{-1}(\D\cup\mo)$. Its complementary $\mathcal{U}_1\backslash N$ is open and
intersects all the $1$-dimensional fibers of 
$\mathcal{U}_1\rightarrow\mathcal{U}_0$, then it
is dense. This implies that $N$ is a proper closed subset, equivalently it has
dimension lower than $r+1$. Moreover the inverse image
$M_i=\pi_i^{-1}(N)$ has dimension dim $M_i<$ dim $M=2r$ because also the 
restricted map $\pi_i:M_i\rightarrow N$ is flat of relative dimension $r-1$.   
Let $\tilde{M}_i$ be the resolution of singularities of $M_i$. It has
the same dimension as $M_i$. Set $\Gamma=\prod_{i=1}^{r}\G{p_i}$, for generic
fixed points $p_1,\ldots,p_r\in Q$ and consider the inverse image of $\Gamma$ 
in $\tilde{M}_i$ via the evaluation map, i.e. the composition:
\begin{displaymath}
ev_i:\tilde{M}_i\twoheadrightarrow M_i\rTo^{ev}\h^r
\end{displaymath} 
\noindent
We apply the Position Lemma to $ev_i$ with the group $\mathcal{A}_0$
 acting on $\h$ (by \ref{rem:pos-0}). 
By \ref{rem:weak-pos}, $ev_i^{-1}(\Gamma)$ has pure dimension equal to 
dim $M_i~-$
cod$(\Gamma\subseteq\h^r)<0$, that is to say it is empty. In particular
$ev^{-1}(\Gamma)\cap M_i=\emptyset$. 
\end{proof}
\noindent
We are ready to give a proof of Theorem \ref{main}. 
We will use induction on the
number of components of the source curve $C$ and we will apply the 
Position Lemma with respect to the action of $\mathcal{A}_0$ on $\h$.
\begin{proof}
\textbf{STEP 1.} 
The subset
$M=\{[C,\mu,x_j]: C\cong\mathbb{P}^1, ~\mu(C)\nsubseteq\D\cup\mo\}$ is open
and smooth in $\overline{M}_{0,r}(\h,\beta)$ and we can consider the 
restriction of
the evaluation map $ev:\mbar{0,r}(\h,(a,b,c))\rightarrow\h$ to it.
Set $\Gamma=\prod_{i=1}^{r}\G{p_i}$.
By the Position Lemma, dim $ev^{-1}(\Gamma)=0$ since $M$ is of the 
expected dimension $2r$.\\
\textbf{STEP 2.} Suppose that $C$ is irreducible and $\mu(C)\subseteq\D$, then 
in $\D$ we have $\mu_*[C]
=\tilde{\beta}=(a/2,b/2,c)$ by what we showed in 
\S\ref{sec:delta}. Let $a'=a/2, b'=b/2$. 
The image of $\mu(C)$ via the support map is a curve $B$ of genus zero and  
bi-degree $(a',b')$ on $Q$.  
The cycles $\Gamma(p_i)$ restricted to $\D$ have codimension 
$2$ and the curve $\mu(C)$ is incident to all of them if and only
if the image curve 
$B$ goes through all the points $p_i$. A rational curve on $Q$ of bi-degree
 $(a',b')$ passes through at most $s$ generical points of $Q$, where:
\begin{displaymath}
2s=\textrm{dim}~Q+\int_{(a',b')}c_1(T_Q)-3+s \Rightarrow s=2a'+2b'-1=a+b-1
\end{displaymath} 
\noindent
We have $s<r=2a+2b+1$, so the irreducible curves $\mu(C)\subseteq \D$ give no 
contribution to our calculations.\\
\textbf{STEP 3.} Now we analyse the contribution from irreducible rational 
curves $C$ such that $\mu(C)\subseteq\mo$. Since $\mo$ is the disjoint union 
$\wt{1}\sqcup\wt{2}$ and $\mu(C)$ is irreducible, it is enough to consider 
the case $\mu(C)\subseteq\wt{1}$. The pushforward class $\mu_*[C]$ in $\wt{1}$
is $(a,b/2)$ with $b$ even. Let $\varphi_1^r$ be the map induced by 
the restricted blowup $\varphi_1:\wt{1}\rightarrow\w{1}$. 
Then we have a composition map:
\begin{displaymath} 
\overline{M}_{0,r}(\wt{1},(a,b/2))\rTo^{ev}\wt{1}^r
\rTo^{\varphi_1^r} W_1^r\subseteq\g^r
\end{displaymath}
\noindent
If a curve of class $(a,b/2)$ intersects all the cycles $\Gamma(p_i)$ then its
image via $\varphi_1$ is of class
$(\varphi_1)_*(a,b/2)=b/2\cdot [W_1]=b[\sig{2,1}]$ because $W_1$ is a
quadric in $\g$, and it goes through all the points $l_1(p_i)\in\g$. Such a
curve passes through at most $s$ fixed points in $\g$, with $s$ given by the
formula:
\begin{displaymath} 
4s=\textrm{dim}~\g+\int_{b[\sig{2,1}]} c_1(T_{\g})-3+s \Rightarrow 
s=\frac{1+4b}{3}
\end{displaymath}
Since $s<r$ we verify that irreducible curves mapped into $\mo$ 
give no contribution to our computation.\\
Suppose that $C$ is the union of $k$ irreducible components and 
$\mu(C)\subseteq\mo$. Since $\wt{1},\wt{2}$ are disjoint, if an irreducible
component is mapped into $\wt{i}$ then all the components are actually mapped
into the same divisor $\wt{i}$, by connectedness. We can assume
$\mu(C)\subseteq\wt{1}$. The number $k$ of components is bounded. In fact
$\mu_*[C]=(a,b,b)$ in $\h$, with $b$ even, then $k$ is at most equal to
$a+b/2$. This implies that $\mu(C)$ goes through at most
$s=\frac{k+4b}{3}\leq \frac{2a+9b}{6}$ cycles. We get $s<r$ also in this case.
\\
Lemma \ref{lem:disjoint} concludes the proof of 1.-2.\\
\textbf{STEP 4.} Suppose $C$ is reducible and $\mu(C)\subseteq\D\cup\mo$. In
particular assume that $C$ has $k$ irreducible components $C_i$ such that:
\begin{displaymath}
\begin{array}{ll}    
C_i\subseteq\D &\textrm{for}~ 1\leq i\leq k_1\\
C_i\subseteq\wt{1} &\textrm{for}~ k_1+1\leq i\leq k_2\\
C_i\subseteq\wt{2} &\textrm{for}~ k_2+1\leq i\leq k
\end{array}
\end{displaymath}
\noindent
We fix the notations:
\begin{displaymath}
\begin{array}{ll}
D_1=\bigcup_{i=1}^{k_1} C_i &\textrm{is of class} ~(a_1,b_1,c_1)\\
 & \\
D_2=\bigcup_{i=k_1+1}^{k_2} C_i &\textrm{is of class} ~(a_2,b_2,c_2)\\ 
 & \\
D_3=\bigcup_{i=k_2+1}^{k} C_i &\textrm{is of class} ~(a_3,b_3,c_3)
\end{array}
\end{displaymath}
\noindent
The conditions $\sum a_j=a,~\sum b_j=b,~\sum c_j=c$ hold. 
The image curve $\mu(D_j)$ intersects $r_j$
cycles. By the previous results we know that $r_j\leq 2a_j+2b_j$ for all $j$
then $r_1+r_2+r_3\leq 2a+2b<r$. The curve $C$ does not intersects all the
cycles $\Gamma(p_i)$.\\
\textbf{STEP 5.} Let $R$ be the subset of $\mbar{0,r}(\h,(a,b,c))$
parametrizing stable maps $[C,\mu,x_j]$ such that $C=\bigcup C_i$ and each
$C_i$ intersects $\mo_4$, the dense orbit. It is a proper closed subset of 
the smooth locus, hence it has dimension lower than $2r$. 
For generic points $p_i$
 the intersection $R\cap ev^{-1}(\Gamma)$ is empty. \\
Finally we
 analyse the contribution from stable maps $\mu:C\rightarrow\h$ with
 rational reducible domain and such that there exists at
 least one component of $C$ mapped into $\D\cup\mo$. We can write $C=C_0\cup
 C_1$ with $C_0\cap C_1=\{p\}$ a point mapped in $\D\cup\mo$.
 Set $[\mu(C_i)]=(a_i,b_i,c_i)$, with 
$\sum a_i=a,~\sum b_i=b,~\sum c_i=c$.\\
We suppose $(a_i,b_i)\neq (0,0)$ for $i=0,1$ and let $\mu(C)$ be incident to
 all the cycles $\Gamma(p_i)$. We know that $\mu(C_i)$ intersects $r_i=
2a_i+2b_i+1-k_i$ cycles, with $k_0+k_1\leq 1$. Assume $r_0=2a_0+2b_0+1$
 and $r_1=2a_1+2b_1$. Theorem \ref{main} applies to $C_0$, by induction.  
It implies that $\mu(C_1)$ intersects $r_1+1$ cycles and a point of 
intersection is in $\D\cup\mo$ (the proof is exactly the same as in \cite{gr} 
Thm. 2.7). This is impossible. Then $a_1=b_1=0$.    
\end{proof}
\begin{rem}
Note that because of smoothness in Corollary \ref{cor:comb} hyperelliptic 
curves are enumerated with multiplicity one. This fact was not proven in 
\cite{p}. 
\end{rem}
\small{
}

\end{document}